\documentclass[final,3p,times,fleqn]{elsarticle}
\usepackage{amssymb} 
\usepackage{amsmath} 
\usepackage{dcolumn}% Align table columns on decimal point
\usepackage{bm}% bold math
\usepackage{booktabs}
\usepackage{amssymb}
\usepackage{amsmath}
\usepackage{array}
\usepackage{colortbl}
\usepackage{color} 
\usepackage{tabularx}
\usepackage{multirow} 
\usepackage{graphicx}
\usepackage{epstopdf}
\usepackage{subfig} 
\usepackage{bbding}
\usepackage{pifont}
\usepackage{wasysym}
\usepackage{utfsym}
\usepackage{fontawesome}

\usepackage{xcolor}
\usepackage{algorithm}
\usepackage{algpseudocode} 
\begin{document}
	
	\begin{frontmatter}
		
		\title{An automatic approach to develop the fourth-order and $L^2$-stable lattice Boltzmann model for diagonal-anisotropic diffusion equations}
		
		\author[a,b,c]{Ying Chen} 
		\author[a,b,c,d]{Zhenhua Chai \corref{cor1}}
		\ead{hustczh@hust.edu.cn}
		\author[a,b,c]{Baochang Shi}
		\address[a]{School of Mathematics and Statistics, Huazhong University of Science and Technology, Wuhan 430074, China}
		\address[b]{Institute of Interdisciplinary Research for Mathematics and Applied Science, Huazhong University of Science and Technology, Wuhan 430074, China}
		\address[c]{Hubei Key Laboratory of Engineering Modeling and
			Scientific Computing, Huazhong University of Science and Technology,
			Wuhan 430074, China}
		\address[d]{The State Key Laboratory of Intelligent Manufacturing Equipment and Technology, Huazhong University of Science and Technology, Wuhan 430074, China }
		\cortext[cor1]{Corresponding author.}
		\begin{abstract} 
	This paper discusses how to develop a high-order multiple-relaxation-time lattice Boltzmann (MRT-LB) model for the general $d$($\geq 1$)-dimensional diagonal-anisotropic diffusion equation. Such an MRT-LB model considers the transformation matrix constructed in a natural way and the D$d$Q($2d^2+1$) [($2d^2+1$) discrete velocities in $d$-dimensional space] lattice structure. A key step in developing the high-order MRT-LB model is to determine the additional adjustable  relaxation parameters  and weight coefficients, which are used to eliminate the truncation errors at some certain orders of the MRT-LB model, while ensuring the stability of  the MRT-LB model.  In this work, we first present a unified MRT-LB model for the $d$-dimensional diagonal-anisotropic diffusion equation. Then, through the direct Taylor expansion, we analyze the macroscopic modified equations of the MRT-LB model up to fourth-order at the diffusive scaling, and further derive the   conditions that ensure  the MRT-LB model to be fourth-order consistent with the diagonal-anisotropic diffusion equation. In particular, when the diagonal-anisotropic diffusion equation is reduced to the isotropic type, we propose another MRT-LB model with the D$d$Q($2d+1$) lattice structure [fewer discrete velocities than the D$d$Q($2d^2+1$) lattice structure], and the fourth-order conditions are similarly derived. Additionally, we also construct the fourth-order initialization scheme for the present LB method. After that,  the condition which guarantees that the MRT-LB model can satisfy the stability structure is explicitly given, and we would like to point out that from a numerical perspective, once the stability structure is satisfied, the MRT-LB model must be $L^2$ stable. In combination with the fourth-order consistent and $L^2$ stability conditions, the relaxation parameters and weight coefficients of the MRT-LB model can be  automatically given  by a simple computer code. Finally, we perform numerical simulations of several benchmark problems, and find that the numerical results can achieve a fourth-order convergence rate, which is in agreement with our theoretical analysis. In particular,  for the isotropic diffusion equation, we also make a comparison between the fourth-order MRT-LB models with the D$d$Q($2d^2+1$) and D$d$Q($2d+1$) lattice structures, and the numerical results show that the MRT-LB model with the D$d$Q($2d^2+1$)  lattice structure is more general.
		\end{abstract}

		
	\end{frontmatter}
	
	\section{Introduction}\label{introduction}
 
It is well known that diffusion equations control numerous physical processes, such as mass transport in membranes \cite{MAVROUDI2006103}, heat transfer in extended surfaces \cite{10.1115/1.3450546} or domains with internal heat sources \cite{Estrada1991}, population evolution in biological systems \cite{CHU20071618}, and the spread and prediction of infectious diseases \cite{Yanke2014PatternFI}, ect. When some specific initial and boundary conditions are imposed on the diffusion equation, the analytical solution can sometimes be obtained \cite{Crank1980}. However, it is more convenient to study it by the numerical methods. In this work, we will  consider the following general $d$($\geq1$)-dimensional diagonal-anisotropic diffusion equation, 
\begin{align}\label{DE}
	\begin{cases}
		\partial_t\phi(\mathbf{x},t)=\sum_{i=1}^d\kappa_{x_i}\partial_{x_i}^2\phi(\mathbf{x},t)+R(\mathbf{x},t),&\mathbf{x} \in\Omega \subset\mathbb{R}^d,\quad t\in (0,T],\\
		\phi(\mathbf{x},0)=\phi^0(\mathbf{x}),&\mathbf{x}\in\Omega,
	\end{cases} 
\end{align}
here, the variable $\phi$ represents the diffusive transport of mass, heat, or any other scalar quantity, $\phi^0$ is the initial value, and both $\phi$ and $\phi^0$ are assumed to be given smooth functions. The constant $\kappa_{x_i}\in\mathbb{R}$ indicates the diffusion coefficient in the $i$-th spatial direction, $R=\eta\phi+S$  with $\eta,S\in\mathbb{R}$ is a linear term. And we note that the physical domain $\Omega$ subjected to the periodic  boundary condition is taken into account in this study.

Over the past few years, some traditional macroscopic numerical methods, such as the finite-difference method \cite{Zhang1998}, the finite volume method \cite{Codina1998},  and the finite element method \cite{Theeraek2011}, have been developed to solve the diffusion equations, where the  discretization process is performed at the macroscopic level. In contrast to these macroscopic numerical methods, the lattice Boltzmann (LB) method, which has been widely used in the last three decades, is an efficient second-order mesoscopic approach based on kinetic theory for the Navier-Stokes equations \cite{Guo2013,Kruger2017,Wang2019}, and in particular it has also been successfully extended to solve partial differential equations containing diffusion equations \cite{HUBER20107956,ANCONA1994107,Suga,Lin2022,SILVA2023105735,Chen2023,VANDERSMAN2000766,GINZBURG20051171,RASIN2005453,PhysRevE.79.016701,Chai2013,Chai2016,Chen2022}. What is more,  in recent years, some high-order LB models have also been proposed for the low-dimensional  convection-diffusion equation and diffusion equation.  For the one-dimensional convection-diffusion equation,  a fourth-order multiple-relaxation-time LB (MRT-LB) model with the D1Q3 lattice structure is developed by Chen et al. \cite{Chen2022}. For the one-dimensional diffusion equation, Suga \cite{Suga} first proposed a fourth-order single-relaxation-time LB (SRT-LB) model with the D1Q3 lattice structure. Then, for the same problem, Lin et al. \cite{Lin2022} further considered the MRT collision operator, and on this basis, the LB method can achieve a sixth-order accuracy. In a recent work, unlike the above two works \cite{Suga,Lin2022}, Silva \cite{SILVA2023105735} focused on the one-dimensional equation with a linear source term. In order to describe this reaction-diffusion system efficiently from a numerical point of view,  compared to the traditional diffusion equation solved in Refs. \cite{Suga,Lin2022}, it is necessary to model the introduced solution-dependent source term more accurately. In order to solve this problem, a fourth-order two-relaxation-time LB (TRT-LB) model is presented in this work \cite{SILVA2023105735}.
 Subsequently, for the two-dimensional diffusion equation, Chen et al. \cite{Chen2023} developed a fourth-order MRT-LB model with the D2Q5 lattice structure and also extended this result to the diffusion equation with a linear source term.  Apart from the aforementioned high-order LB models, which are limited to the one- and two-dimensional problems. Recently, Chen et al. developed a unified fourth-order MRT-LB model for the $d$-dimensional coupled Burgers' equations \cite{chen2023cole} and also a unified fourth-order SRT-LB model for the $d$-dimensional hyperbolic equation \cite{chen2024hy}. Despite the valuable work listed above, it can be seen that for the more general $d$-dimensional diagonal-anisotropic diffusion equation with a linear source term (\ref{DE}), it is still unclear how to develop a high-order LB method. Nevertheless, according to the previous works \cite{ANCONA1994107,Suga,Lin2022}, it can be observed that the choice of collision operator and lattice structure would have a significant impact on improving the accuracy of LB methods, due to the introduction of some additional adjustable parameters (e.g, the relaxation parameters in the MRT collision operator \cite{dHumières2002,Chai2020,Chai2023} and the weight coefficients in the lattice structure \cite{Kruger2017}). Therefore, in this study, to develop a higher-order LB model for (\ref{DE}), we will adopt the more general MRT collision operator constructed in a natural way  and the D$d$Q($2d^2+1$) lattice structure \cite{Chai2020,chen2023cole,chen2024hy}.

When developing the LB method with enhanced accuracy requirements, it is essential not only to focus on achieving a high-order accuracy, but also to consider the stability of the LB method, especially for the more widely used MRT-LB model. In Ref. \cite{PhysRevE.61.6546}, the stability of the D2Q9 orthogonal MRT-LB model has been extensively studied using the von Neumann stability theory.  The core idea of the von Neumann stability theory is to derive the growth matrix, which can be obtained by linearizing the LB method at a specific velocity. However, it is typically difficult to specify the conditions under which the moduli of the eigenvalues of the growth matrix are not greater than 1, i.e., the explicit stability conditions for the LB method. So far, the rigorous studies on the stability analysis and the determination of its explicit conditions are scarcely available, with the only research works being in Refs. \cite{Bellotti2022,Chen2022,Caetano2023,Chen2023}. In general, the stability conditions, which are associated with the relaxation parameters of the MRT collision operator, are usually obtained by using a computer code.   Apart from the stability analysis based on the von Neumann stability theory, a weighted $L^2$-stability result for the LB method linearized at zero velocity has been recently validated in Ref.  \cite{JunkYong2009}. And it is worth noting that in terms of the weighted $L^2$ stability analysis, not only the stability conditions can be explicitly obtained from validating whether the LB method satisfies the stability structure \cite{Rheinländer2010,Yong2009}, but also the convergence proof of the nonlinear problems can be established \cite{junk2008convergence,JunkYang2009}.  Nevertheless, such stability has only been demonstrated for the D2Q9 orthogonal MRT model and how to prove it for many other MRT models is still unclear.  To bridge this knowledge gap, Yong et al. \cite{Yang2024} put forward an automatic approach for the stability analysis of MRT-LB models and applied it to ten different MRT-LB models \cite{dHumières1992,dHumières2002,Suga2015,Kaehler2013,Liu2016,DeRosis2020,Fakhari2017}. The key step in this approach is to decompose the Jacobian matrix of the collision term in the LB model into the product of a symmetric semi-negative definite matrix and a diagonally positive definite matrix, in particular, this decomposition process can be automatically verified by a simple computer code. With the aid of this work \cite{Yang2024}, in our current study, we will develop a high-order LB model for   (\ref{DE}) and present its explicit stability structure preserving condition.

	The remainder of this paper is structured as follows. In Section \ref{MRT-LB-DE}, we present the MRT-LB model for the $d$-dimensional diagonal-anisotropic diffusion equation. In Section \ref{Fourth-MRT-LB}, through the direct Taylor expansion, we first deduce the conditions that ensure the present MRT-LB model to be fourth-order consistent with the  diagonal-anisotropic diffusion equation, and then present the fourth-order   initialization scheme for the MRT-LB model. Thereafter, the condition which guarantees that the MRT-LB model can satisfy the stability structure is provided, and the relationship between the stability structure preserving condition and the $L^2$ stability of the MRT-LB model is also discussed. Some numerical experiments are carried out in Section \ref{Numer}, and finally some conclusions are summarized in Section \ref{conclusion}.

	\section{The MRT-LB model for the diagonal-anisotropic diffusion equation}\label{MRT-LB-DE} 
	In this section, we will present a unified MRT-LB model for the diagonal-anisotropic diffusion equation (\ref{DE}), where the transformation matrix constructed in a natural way and the D$d$Q$(2d^2+1)$ lattice structure are employed. In particular, for the isotropic diffusion equation, another MRT-LB model with the D$d$Q$(2d+1)$ lattice structure [fewer discrete velocities than the D$d$Q$(2d^2+1)$ lattice structure] is also provided.

	\subsection{Spatial and temporal discretization}
	For the sake of brevity and to facilitate the following analysis, the physical domain $\Omega$ is first discretized into a uniform mesh $\mathcal{L}$ with a lattice spacing $\Delta x>0$, where $\mathcal{L}=\Delta x\mathbb{Z}^d$, and the corresponding lattice node is denoted by $\mathbf{x}_i$. The time is uniformly discretized as $t_n=n\Delta t$, where $\Delta t>0$ represents the time step. The lattice velocity in the LB method can then be expressed as $c = \Delta x/\Delta t$. For the diagonal-anisotropic diffusion equation (\ref{DE}), the diffusive scaling ($\Delta t\propto \Delta x^2$) is adopted in this work. Without loss of generality, we can also take $\Delta t=\xi\Delta x^2$, where the parameter $\xi\in\mathbb{R}^+$ and $\xi=0$ when $\Delta x$, $\Delta t\rightarrow 0$.
	\subsection{The MRT-LB model  }\label{MRT-LB}
	To develop a high-order LB method for the diagonal-anisotropic diffusion equation (\ref{DE}), we here adopt the general MRT-LB model  \cite{dHumières2002,Chai2020,Chai2023} with some additional adjustable relaxation parameters, which can be used to eliminate some certain high-order truncation errors of the LB method,  in which the evolution equation can be written as
	\begin{align}\label{lb}
		f_k(\mathbf{x}+\mathbf{c}_k\Delta t,t+\Delta t)=f_k(\mathbf{x},t)-\sum_{i=1}^q\bigg[\bm{\Lambda}_{ki}\big(f_i-f_i^{eq}\big)\bigg](\mathbf{x},t)+\Delta t\sum_{i=1}^q\bigg[\Big(\mathbf{I}-\frac{\bm{\Lambda}}{2}\Big)_{ki}R_i\bigg](\mathbf{x},t),\quad k\in \llbracket1,q\rrbracket,\footnotemark[1]
	\end{align}  where $f_k(\mathbf{x},t)$, $f_k^{eq}(\mathbf{x},t)$, and $R_k$  represent the distribution function, equilibrium distribution function, and discrete source term at position $\mathbf{x}$ and time $t$, respectively. In terms of the present MRT-LB model (\ref{lb}) for  (\ref{DE}), we consider the D$d$Q$q$ lattice structure with the number of the discrete velocities $q=(2d^2+1)$. To be specific,  the velocity set $\mathbf{c}$ of the D$d$Q$q$ lattice structure can be given by
		\begin{subequations}\label{DdQq}
		\begin{align}    
				d=1:\quad&\mathbf{c}=\left(\begin{matrix}
					0&1&-1\end{matrix}\right)^Tc, \\
			 	d=2:\quad&		  
				\mathbf{c} =\left(\begin{matrix}
					0&1&0&-1&0&1&-1&-1&1\\
					0&0&1&0&-1&1&1&-1&-1
				\end{matrix}\right)^Tc,\\
		 	d=3:\quad&\mathbf{c} =\left(\begin{array}{ccccccccccccccccccc}
					0&1&0&0&-1&0&0&1&-1&-1& 1&\\ 
					0&0&1&0&0&-1&0&1& 1&-1&-1&\\
					0&0&0&1&0&0&-1&0& 0& 0& 0&
				\end{array}\right.\notag\\
				&\qquad\qquad  \left.\begin{array}{ccccccccccccccccccc}
					1&-1&-1& 1&0&0&0&0\\ 
					0& 0& 0& 0&1&-1&-1&1\\
					1& 1&-1&-1&1& 1&-1&-1
				\end{array}\right)^Tc, 
		\end{align}
	\end{subequations}
	and the discrete velocity $(\mathbf{c}_k)^T$ in (\ref{lb}) is then given by the $k$-th row of the above velocity set $\mathbf{c}$. The collision matrix $\bm{\Lambda}=\mathbf{M}^{-1}\mathbf{SM}\in\mathbb{R}^{q\times q}$ in (\ref{lb}) is an invertible matrix, where the transformation matrix $\mathbf{M}$ constructed in a natural way and the   relaxation matrix $\mathbf{S}$ are respectively given by
 	\footnotetext[1]{The notation $\llbracket{1,q}\rrbracket$ denotes $\{1,2,\cdots,q\}$.}
	\begin{subequations}\label{DdQq-M}
		\begin{flalign}
			&&d=1:\:\:   
			&\begin{cases}\label{d1q3} 
				\mathbf{p}_{\mathbf{M}}(X_1)=\big(1,X_1,X_1^2\big),  \\
				\mathbf{S}=\makebox{\textbf{diag}}\big(s_0,s_{x_1},s_{2|x_1^2}\big), 
			\end{cases} &\\
			&&	d=2:\:\:  				&\begin{cases}\label{d2q9}
				\mathbf{p}_{\mathbf{M}}(X_1,X_2)=\Big(1,X_1,X_2,X_1^2,X_2^2,X_1X_2,X_1^2X_2,X_1X_2^2,X_1^2X_2^2\Big),   \\
				\mathbf{S}=\makebox{\textbf{diag}}\big(s_0,s_{x_1},s_{x_2},s_{2|x_1^2},s_{2|x_2^2},s_{2|x_1x_2},s_{3|x_1^2x_2},s_{3|x_1x_2^2},s_{4|x_1^2x_2^2}\big), 
			\end{cases} &\\
			&&	d=3:\:\:	
				&\begin{cases}\label{d3q19}
					\mathbf{p}_{\mathbf{M}}(X_1,X_2,X_3)=\Big(1,X_1,X_2,X_3,X_1^2,X_2^2,X_3^2,X_1X_2,X_1X_3,X_2X_3,\\
				\qquad\qquad\qquad\qquad\quad X_1^2X_2,X_1^2X_3,X_2^2X_1,X_2^2X_3,
				X_3^2X_1,X_3^2X_2,\\
				\qquad\qquad\qquad\qquad\quad X_1^2X_2^2,X_1^2X_3^2,X_2^2X_3^2\Big) , \\
				\mathbf{S}=\makebox{\textbf{diag}}\big(s_0,s_{x_1},s_{x_2},s_{x_3},s_{2|x_1^2},s_{2|x_2^2},s_{2|x_3^2},s_{2|x_1x_2},s_{2|x_1x_3},s_{2|x_2x_3},s_{3|x_1^2x_2},s_{3|x_1^2x_3}, \\
				\qquad\qquad\qquad  s_{3|x_1x_2^2},s_{3|x_2^2x_3},s_{3|x_1x_3^2},s_{3|x_2x_3^2},s_{4|x_1^2x_2^2},s_{4|x_2^2x_3^2}\big), 
			\end{cases} &
		\end{flalign}
	\end{subequations}
	here $\mathbf{p}_{\mathbf{M}}(X_1,X_2,\cdots,X_d)$ is a polynomial set  of the transformation matrix $\mathbf{M}$, and the $i$-th ($i\in\llbracket1,q\rrbracket$) row of  $\mathbf{M}$ is determined by the $i$-th element of   $\mathbf{p}_{\mathbf{M}}$ with $X_k=(\mathbf{c}^T)_k$ ($k\in\llbracket 1,d\rrbracket$). The  velocity set $\mathbf{c}$, the transformation matrix $\mathbf{M}$, and the diagonal relaxation matrix $\mathbf{S}$ for the case $d>3$   can be found in our previous works \cite{chen2023cole,chen2024hy}.  
	
To correctly recover (\ref{DE}) from the present MRT-LB model (\ref{lb}),  the equilibrium distribution function  and  discrete source term  in (\ref{lb}) can be designed as   
			\begin{align*}
f_k^{eq}=\omega_{k-1}\phi,\quad R_k=\omega_{k-1}R,  
			\end{align*} 
		where $\omega_k\in(0,1)$ for all $k\in\llbracket 0,q-1\rrbracket$ are the weight coefficients, and in our work, they are assumed to satisfy the following relations
		\begin{subequations}\label{weight}
			\begin{align}
				&\omega_i=\omega_{i+d},\quad i\in\llbracket 1,d\rrbracket\\
				&\omega_{2d+1}=\omega_{2d+2}=\cdots=\omega_{2d^2}\overset{\Delta}{=}\tilde{\omega},\\
				&\omega_0=1-2\sum_{i=1}^{d}\omega_i-2d(d-1)\tilde{\omega}. 
			\end{align}
		\end{subequations}

Then, through some commonly used asymptotic methods \cite{chapman1990,HOLDYCH2004,d2009viscosity,Ginzburg2012,Dubois2022}, it is easy to show that in the sense of second-order accuracy,  the macroscopic variable $\phi$ of (\ref{DE}) can be calculated from the present MRT-LB model (\ref{lb}) by 
\begin{align*}
	\phi=\frac{2\sum_{k=1}^qf_k+\Delta tS}{2-\Delta t\eta},
\end{align*}
where the fact of $\sum_{k=1}^qf_k=\phi-\Delta t(\eta\phi+S)/2$ has been used, and  (\ref{DE})  can be recovered   as long as the diffusion coefficient satisfies 
\begin{align}\label{kappa}
\kappa_{x_i}=\big(2\omega_i+4(d-1)\tilde{\omega}\big)\Big(\frac{1}{s_{x_i}}-\frac{1}{2}\Big)\frac{\Delta x^2}{\Delta t},\quad i\in\llbracket 1,d\rrbracket,
\end{align}
where the relaxation parameter $s_{x_i}$ is located in the range of $(0,2)$.

For the present MRT-LB model (\ref{lb}), we now give two remarks.\\
\textbf{Remark 1.} Similar to our previous works \cite{Chen2022,Chen2023}, one can conclude that the relaxation parameter $s_0$, which corresponds to the zeroth (or conservative) moment $\phi$, can be chosen arbitrarily except zero. And from (\ref{kappa}), which establishes the relationship among $\kappa_{x_i}$, $s_{x_i}$ (related to the first-order moment of the distribution functions), $\omega_i$, and $\tilde{\omega}$, it can be found that apart from $s_{x_i}$, $\omega_i$, and $\tilde{\omega}$, the relaxation parameters in the diagonal relaxation matrix $\mathbf{S}$ which correspond to the second- to fourth-order moments remain unrestricted. This means that these free relaxation parameters could be used to improve the numerical accuracy and/or stability of the MRT-LB model (see the following Parts \ref{accuracy} and \ref{stability} for details). \\
\textbf{Remark 2.} The lattice structure of the LB method is not unique \cite{Kruger2017}. However, we would like to point out that to develop a high-order LB method for the diagonal-anisotropic diffusion equation (\ref{DE}), the MRT-LB model with another commonly used D$d$Q($2d+1$) lattice structure is not appropriate unless (\ref{DE}) is reduced to the isotropic type, and the reason will be explained in the following Part \ref{accuracy}. To be specific,  the velocity set $\mathbf{c}$, the transformation matrix $\mathbf{M}$, and the relaxation matrix $\mathbf{S}$ of the MRT-LB model with the D$d$Q($2d+1$) lattice structure  can be given by
\begin{subequations} 
	\begin{flalign}
		&&d=1:\:\:  &
		\begin{cases} 
			\mathbf{c}=\big(0\quad 1\quad-1\big)^Tc, \\ 
			\mathbf{p}_{\mathbf{M}}(X_1)=\big(1,X_1,X_1^2\big),  \\
			\mathbf{S}=\makebox{\textbf{diag}}\big(s_0,s_{x_1},s_{2|x_1^2}\big), 
		\end{cases} &\\
		&&	d=2:\:\:  &\begin{cases} 
			\mathbf{c} =\left(\begin{matrix}
				0&1&0&-1&0&1&-1&-1&1\\
				0&0&1&0&-1&1&1&-1&-1
			\end{matrix}\right)^Tc,\\ 
			\mathbf{p}_{\mathbf{M}}(X_1,X_2)=\Big(1,X_1,X_2,X_1^2,X_2^2\Big),  \\
			\mathbf{S}=\makebox{\textbf{diag}}\big(s_0,s_{x_1},s_{x_2},s_{2|x_1^2},s_{2|x_2^2},s_{2|x_1x_2}\big), 
		\end{cases} &\\
		&&	d=3:\:\:	&\begin{cases} 
			\mathbf{c} =\left(\begin{array}{ccccccccccccccccccc}
				0&1&0&0&-1&0&0\\ 
				0&0&1&0&0&-1&0\\
				0&0&0&1&0&0&-1
			\end{array}\right).\\ 
			\mathbf{p}_{\mathbf{M}}(X_1,X_2,X_3)=\Big(1,X_1,X_2,X_3,X_1^2,X_2^2,X_3^2\Big) ,\\
			\mathbf{S}=\makebox{\textbf{diag}}\big(s_0,s_{x_1},s_{x_2},s_{x_3},s_{2|x_1^2},s_{2|x_2^2},s_{2|x_3^2},s_{2|x_1x_2},s_{2|x_1x_3},s_{2|x_2x_3}\big)
		\end{cases} &\\
		&&		d>3:\:\:&\begin{cases}
			\mathbf{c} =\left(\begin{matrix}
			\mathbf{0}_{d\times 1}&\mathbf{I}_{d\times d}&-\mathbf{I}_{d\times d}\\
		\end{matrix}\right)^Tc,    \\
		\mathbf{p}_{\mathbf{M}}(X_1,X_2,\cdots,X_d) =\Big(1,X_1,X_2,\cdots,X_d,X_1^2,X_2^2,\cdots,X_d^2 \Big), \\
		\mathbf{S}=\makebox{\textbf{diag}}\big(s_0,s_{x_1},s_{x_2},\cdots,s_{x_d},s_{2|x_1^2},s_{2|x_2^2},\cdots,s_{2|x_d^2},s_{2|x_1x_2},s_{2|x_1x_3},\cdots,s_{2|x_{d-1}x_d}\big).
	\end{cases} &
\end{flalign}
	\end{subequations}

	\section{The fourth-order MRT-LB model for the diagonal-anisotropic diffusion equation}\label{Fourth-MRT-LB} 
		In this section, an accuracy analysis on the MRT-LB model (\ref{lb}) using the direct Taylor expansion is first carried out, and then  the conditions that ensure    the   MRT-LB model (\ref{lb}) to be fourth-order consistent with (\ref{DE}) will be presented. Subsequently, to prevent the reduction of the accuracy order of the MRT-LB model (\ref{lb}), a fourth-order  initialization scheme will also be provided. After that, with the aid of the previous works related to the  stability structure \cite{Banda2006,JunkYang2009,Yang2024}, the stability of the MRT-LB model (\ref{lb}) will be analyzed. Finally, based on the above research, an automatic approach to determine the relaxation parameters and weight coefficients  of the MRT-LB model (\ref{lb}) with a fourth-order accuracy while ensuring its stability will be proposed.
	\subsection{The accuracy analysis of the MRT-LB model}\label{accuracy}
      	
   	Under the premise of $\Delta t=\xi\Delta x^2$,  we now expand the MRT-LB model (\ref{lb}) up to the order of $O(\Delta x^6)$ at position $\mathbf{x}$ and time $t$ [to simplify the following analysis, we denote $\psi=\psi(\mathbf{x},t)$ ($\phi=f_k,f_k^{eq}$, and $R_k$)],
	\begin{align}\label{f-o6}
		\xi\Delta x^2&\partial_tf_k+\frac{\xi^2\Delta x^4}{2}\partial_t^2f_k+O(\Delta x^6)\notag\\
		=&-\Delta x\mathbf{e}_k\cdot\nabla f_k^{eq}+\frac{\Delta x^2}{2}\mathbf{e}_k^2\overset{2}{\cdot}\nabla^2 f_k^{eq}-\frac{\Delta x^3}{6}\mathbf{e}_k^3\overset{3}{\cdot}\nabla^3 f_k^{eq}
		+\frac{\Delta x^4}{24}\mathbf{e}_k^4\overset{4}{\cdot}\nabla^4 f_k^{eq}
		-\frac{\Delta x^5}{120}\mathbf{e}_k^5\overset{5}{\cdot}\nabla^5 f_k^{eq}\notag\\
		&-\sum_{i=1}^q\bm{\Lambda}_{ki}f_i^{ne}+\Delta x\mathbf{e}_k\cdot\nabla \sum_{i=1}^q(\bm{\Lambda}-\mathbf{I})_{ki}f_i^{ne}
		-\frac{\Delta x^2}{2}\mathbf{e}_k^2\overset{2}{\cdot}\nabla^2 \sum_{i=1}^q(\bm{\Lambda}-\mathbf{I})_{ki}f_i^{ne}\notag\\
		&
		+\frac{\Delta x^3}{6}\mathbf{e}_k^3\overset{3}{\cdot}\nabla^3 \sum_{i=1}^q(\bm{\Lambda}-\mathbf{I})_{ki}f_i^{ne}
		-\frac{\Delta x^4}{24}\mathbf{e}_k^4\overset{4}{\cdot}\nabla^4 \sum_{i=1}^q(\bm{\Lambda}-\mathbf{I})_{ki}f_i^{ne}+\xi\Delta x^2\sum_{i=1}^q\Big(\mathbf{I}-\frac{\mathbf{\Lambda}}{2}\Big)_{ki}R_i\notag\\
		& 
		-\xi\Delta x^3\mathbf{e}_k\cdot\nabla\sum_{i=1}^q\Big(\mathbf{I}-\frac{\mathbf{\Lambda}}{2}\Big)_{ki}R_i
		+\frac{\xi\Delta x^4}{2}\mathbf{e}_k^2\overset{2}{\cdot}\nabla^2\sum_{i=1}^q\Big(\mathbf{I}-\frac{\mathbf{\Lambda}}{2}\Big)_{ki}R_i
		-\frac{\xi\Delta x^5}{6}\mathbf{e}_k^3\overset{3}{\cdot}\nabla^3\sum_{i=1}^q\Big(\mathbf{I}-\frac{\mathbf{\Lambda}}{2}\Big)_{ki}R_i,
	\end{align} 
	from which the relation of $f^{ne}_k=f_k-f_k^{eq}=O(\Delta x)$ can be obtained. The velocity $\mathbf{e}_i=(\mathbf{c}_i)^T/c$ and the notation $\mathbf{e}_k^i\overset{i}{\cdot} \nabla^i$ is defined by
	\begin{align*}
		\mathbf{e}_{k_1}^{n_1}\mathbf{e}_{k_2}^{n_2}\cdots \mathbf{e}_{k_m}^{n_m}\overset{i}{\cdot} \nabla^i =\sum_{\alpha_i=1}^{d}\cdots\sum_{\alpha_1=1}^{d}\Bigg(\bigg(\prod_{l=1}^me_{k_l\alpha_1}^{n_l}\bigg)\bigg(\prod_{l=1}^me_{k_l\alpha_2}^{n_l}\bigg)\cdots \bigg(\prod_{l=1}^me_{k_l\alpha_d}^{n_l}\bigg)\nabla^i_{\alpha_1\ldots\alpha_i}\Bigg)  \makebox{   with  } \sum_{l=1}^mn_l=i\in \llbracket2,5\rrbracket.
	\end{align*}
	The key step in deriving the macroscopic modified equation of the mesoscopic MRT-LB model (\ref{lb})  is to calculate the moment of the distribution function $f_k$. Here,  we introduce the   notation $\bm{\varepsilon}_{N_r|x_1^{n_1}x_2^{n_2}\cdots x_d^{n_d}}$  to facilitate the following derivation, 
\begin{align}\label{varepsilon}
	\bm{\varepsilon}_{N_r|x_1^{n_1}x_2^{n_2}\cdots x_d^{n_d}}=\sum_{k=1}^q\bigg(\prod_{l=1}^d\mathbf{c}_{kl}^{n_l}\bigg)\omega_{k-1} \makebox{   with  } N_r=\sum_{l=1}^dn_l,\quad n_l\in\mathbb{R}^+,
\end{align}
which corresponds the $N_r$-th moment of the weight coefficient $\omega_k$. 
On this basis,   the moments of the equilibrium distribution function and discrete source term can be written as
\begin{subequations}\label{meq-R}
	\begin{align}
		&\mathbf{m}^{eq}=\mathbf{Mf}^{eq}=\Big(1,\mathbf{0}_{d\times 1}, \bm{\varepsilon}_{2|\mathbf{x}^{.2}},\mathbf{0}_{d(d-1)/2},\mathbf{0}_{d(d-1)}, \bm{\varepsilon}_{4|(\mathbf{x}^{.2}\mathbf{x}^{.2})_{\alpha<\beta}}\Big)^T  \phi \makebox{   with  } \mathbf{f}^{eq}=\big(f_1^{eq},f_2^{eq},\cdots,f_q^{eq}\big)^T, \\
		&\mathbf{m}^R=\mathbf{MR}=\Big(1,\mathbf{0}_{d\times 1}, \bm{\varepsilon}_{2|\mathbf{x}^{.2}},\mathbf{0}_{d(d-1)/2},\mathbf{0}_{d(d-1)}, \bm{\varepsilon}_{4|(\mathbf{x}^{.2}\mathbf{x}^{.2})_{\alpha<\beta}}\Big)^T R \makebox{   with  } \mathbf{R}=\big(R_1,R_2,\cdots,R_q\big)^T,  
	\end{align}
\end{subequations} 
where the relations of the weight coefficients in (\ref{weight}) have been used, and the two row vectors $\bm{\varepsilon}_{2|\mathbf{x}^{.2}}$ and $\bm{\varepsilon}_{4|(\mathbf{x}^{.2}\mathbf{x}^{.2})_{\alpha<\beta}}$ are given by
\begin{align*} 
	&\bm{\varepsilon}_{2|\mathbf{x}^{.2}} =\Big(\bm{\varepsilon}_{2|x_1^2},\bm{\varepsilon}_{2|x_2^2},\cdots,\bm{\varepsilon}_{2|x_d^2}\Big)\in\mathbb{R}^{1\times d},  \\
	&\bm{\varepsilon}_{4|(\mathbf{x}^{.2}\mathbf{x}^{.2})_{\alpha<\beta}} =\Big(\bm{\varepsilon}_{4|x_1^2x_2^2},\cdots,\bm{\varepsilon}_{4|x_1^2x_d^2},\bm{\varepsilon}_{4|x_2^2x_3^2},\cdots,\bm{\varepsilon}_{4|x_2^2x_d^2},\cdots, \bm{\varepsilon}_{4|x_i^2x_{i+1}^2},\cdots,\bm{\varepsilon}_{4|x_i^2x_d^2},\cdots,\bm{\varepsilon}_{4|x_{d-1}^2x_d^2}\Big)\in\mathbb{R}^{1\times \big(d(d-1)/2\big)}, 
\end{align*}
with
\begin{subequations}\label{moment}
	\begin{align}
		&\bm{\varepsilon}_{2|x_i^2}=\big(2\omega_i+4(d-1)\tilde{\omega}\big)\phi,\quad i\in\llbracket1,d\rrbracket,\\
		&\bm{\varepsilon}_{4|x_i^2x_j^2}=\begin{cases}
			4\tilde{\omega}\phi,& i,j\in\llbracket1,d\rrbracket,\quad i\neq j,\\
			\big(2\omega_i+4(d-1)\tilde{\omega}\big)\phi,& i,j\in\llbracket1,d\rrbracket,\quad i=j.\\
		\end{cases}
	\end{align} 
\end{subequations}
Based on (\ref{weight}), (\ref{varepsilon}), and (\ref{meq-R}), the zeroth moment of (\ref{f-o6}) is expressed as 
	\begin{align} \label{macro-1}
		&\xi\partial_t\phi-\frac{\xi^2\eta\Delta x^2}{2}\partial_t\phi+\frac{\xi^2\Delta x^2}{2}\partial_t^2\phi+O(\Delta x^4)=\frac{1}{2}\sum_{i,j=1}^d\nabla^2_{ij}\bm{\varepsilon}_{2|x_ix_j}\phi
		+\frac{\Delta x^2}{24}\sum_{i,j,l,p=1}^d\nabla^4_{ijlp}\bm{\varepsilon}_{4|x_ix_jx_lx_p}\phi\notag\\
	 &+\frac{1}{\Delta x}\sum_{i=1}^d\nabla_i(s_{x_i}-1)\sum_{k=1}^q\mathbf{e}_{ki}  f_k^{ne}
		-\frac{1}{2}\sum_{i,j=1}^d\nabla^2_{ij}(s_{2|x_ix_j}-1)\sum_{k=1}^q\mathbf{e}_{ki}\mathbf{e}_{kj}f_k^{ne}	+\xi R\notag\\
		&+\frac{\xi\eta\Delta x^2}{2}\sum_{i,j=1}^d\nabla^2_{ij} \Big(1-\frac{s_{2|x_ix_j}}{2}\Big)\bm{\varepsilon}_{2|x_ix_j}\phi+\frac{\Delta x}{6}\sum_{i,j,l=1}^d\nabla^3_{ijl}(s_{3|x_ix_jx_l}-1)\sum_{k=1}^q\mathbf{e}_{ki}\mathbf{e}_{kj}\mathbf{e}_{kl} f_k^{ne}\notag\\
		&
		-\frac{\Delta x^2}{24}\sum_{i,j,l,p=1}^d\nabla^4_{ijlp}(s_{4|x_ix_jx_lx_p}-1)\sum_{k=1}^q\mathbf{e}_{ki}\mathbf{e}_{kj}\mathbf{e}_{kl}\mathbf{e}_{kp}f_k^{ne},
	\end{align} 
where from the relation of $\mathbf{M}\bm{\Lambda}=\mathbf{SM}$, (\ref{DdQq}), (\ref{DdQq-M}), and (\ref{weight}), the parameters  $s_{N_r|x_1^{n_1}x_2^{n_2}\cdots x_d^{n_d}}$ for $N_r=3$ and 4 introduced in the above equation   satisfy 
\begin{align*}
\sum_{k=1}^q\mathbf{e}_{k1}^{n_1}\mathbf{e}_{k2}^{n_2}\cdots \mathbf{e}_{kd}^{n_d}\bm{\Lambda}_{ki}=s_{N_r|x_1^{n_1}x_2^{n_2}\cdots x_d^{n_d}}\mathbf{e}_{i1}^{n_1}\mathbf{e}_{i2}^{n_2}\cdots \mathbf{e}_{id}^{n_d} {\makebox{ with }} \sum_{l=1}^dn_l=N_r,
\end{align*} 
 and they are defined by
\begin{subequations}\label{la-m-s}
	\begin{align}
		&s_{3|x_kx_lx_p}= \begin{cases}
			 s_{x_k},& k\in\llbracket 1, d\rrbracket,\quad p=l=k,\\
			s_{3|x_k^2x_l},& k,l\in\llbracket 1,d\rrbracket,\quad p=k,k<l,\\
			1,& k,l,p\in\llbracket 1,d\rrbracket,\quad k\neq l\neq p,
		\end{cases}\\ 
		& s_{4|x_kx_lx_px_r}= \begin{cases} 
			s_{2|x_k^2},& k\in\llbracket 1,d\rrbracket,\quad r=p=l=k,\\
			s_{2|x_kx_l},& k,l\in\llbracket 1,d\rrbracket,\quad r=p=k,k<l,\\ 
			s_{4|x_k^2x_l^2},& k,p\in\llbracket 1,d\rrbracket,\quad r= k,p=l,k<l,\\ 
			1,& k,l,p\in\llbracket 1,d\rrbracket,\quad r=k,k\neq l\neq p,\\ 
			1,& k,l,p,r\in\llbracket 1,d\rrbracket,\quad k\neq l\neq p\neq r.
		\end{cases} 
	\end{align}
\end{subequations}

It is obvious that from (\ref{la-m-s}), the following four   terms,   
	\begin{align*}
		&\frac{1}{\Delta x}\sum_{k=1}^q\mathbf{e}_{ki}  f_k^{ne}	,\quad\sum_{k=1}^q\mathbf{e}_{ki}\mathbf{e}_{kj}f_k^{ne},\quad \Delta x\sum_{k=1}^q\mathbf{e}_{ki}\mathbf{e}_{kj}\mathbf{e}_{kl} f_k^{ne},\quad\Delta x^2\sum_{k=1}^q\mathbf{e}_{ki}\mathbf{e}_{kj}\mathbf{e}_{kl}\mathbf{e}_{kp}f_k^{ne},
	\end{align*} 
	are unknown, which must be determined with at least fifth-, fourth-, third-, and second-order accuracy, respectively. This indicates that  the second- to fifth-order expansions of the distribution function $f_k$ should be derived. For the sake of brevity, we here only show the final results and the details can be found in \ref{sec-app-1}, 
	\begin{subequations}\label{fneq-moment}
		\begin{align} 
			&  \sum_{k=1}^q\mathbf{e}_{ki}f_k^{ne} = -\frac{\xi\Delta x^3}{s_{x_i}}\partial_t\sum_{j=1}^d\nabla_{j}\Big(1-\frac{1}{s_{2|x_ix_j}}+\frac{1}{s_{x_i}}\Big)\bm{\varepsilon}_{2|x_ix_j}\phi-\frac{\Delta x}{s_{x_i}}\sum_{j=1}^d\nabla_{j}\bm{\varepsilon}_{2|x_ix_j}-\frac{\Delta x^3}{6s_{x_i}}\sum_{j,l,p=1}^d\nabla^3_{jlp}\bm{\varepsilon}_{4|x_ix_jx_lx_p}\phi \notag\\
			&\qquad\qquad\quad-\frac{\Delta x^3}{s_{x_i}}\sum_{j,l,p=1}^d\nabla^3_{jlp}\Big(\frac{1}{s_{3|x_ix_jx_l}s_{2|x_ix_j}}-\frac{1}{2s_{3|x_ix_jx_l}}-\frac{1}{2s_{2|x_ix_j}}\Big) \bm{\varepsilon}_{4|x_ix_jx_lx_p}\phi\notag\\
			& \qquad\qquad \quad
			-\frac{\xi\Delta x^3}{s_{x_i}}\sum_{j=1}^d\nabla_{j}\Big(\frac{1}{s_{2|x_ix_j}}-\frac{1}{2}\Big) \bm{\varepsilon}_{2|x_ix_j}R+O(\Delta x^5),\label{fneq-5rd}\\
			&\sum_{k=1}^q\mathbf{e}_{ki}\mathbf{e}_{kj}f_k^{ne}=-\xi \frac{\Delta x^2}{s_{2|x_ix_j}}\partial_t\bm{\varepsilon}_{2|x_ix_j}\phi+\frac{\Delta x^2}{s_{2|{x_ix_j}}} \sum_{l,p=1}^d\nabla^2_{lp}\Big(\frac{1}{s_{3|x_ix_jx_l}}-\frac{1}{2}\Big)\bm{\varepsilon}_{4|x_ix_jx_lx_p}\phi+\xi\Delta x^2\Big(\frac{1}{s_{2|x_ix_j}}-\frac{1}{2}\Big)\bm{\varepsilon}_{2|x_ix_j}R+O(\Delta x^4)  ,\label{fneq-4rd}\\
			&\sum_{k=1}^q\mathbf{e}_{ki}\mathbf{e}_{kj}\mathbf{e}_{kl}f_k^{ne}=-\frac{\Delta x}{s_{3|x_ix_jx_l}}\sum_{p=1}^d\nabla_{p}\bm{\varepsilon}_{4|x_ix_jx_lx_p}\phi+O(\Delta x^3)  ,\\
			&   \sum_k\mathbf{e}_{ki}\mathbf{e}_{kj}\mathbf{e}_{kl}\mathbf{e}_{kp}f_k^{ne}= O(\Delta x^2)   .   
		\end{align} 
	\end{subequations} 
	Substituting (\ref{fneq-moment}) into (\ref{macro-1}), one can obtain the fourth-order modified equation of the MRT-LB model (\ref{lb}),
	 	\begin{align}\label{err}
	 	&\xi\partial_t\phi-\frac{\xi^2\eta\Delta x^2}{2}\partial_t\phi+\frac{\xi^2\Delta x^2}{2}\partial_t^2\phi+O(\Delta x^4)=\frac{1}{2}\sum_{i,j=1}^d\nabla^2_{ij}\bm{\varepsilon}_{2|x_ix_j}\phi
	 	+\frac{\Delta x^2}{24}\sum_{i,j,l,p=1}^d\nabla^4_{ijlp}\bm{\varepsilon}_{4|x_ix_jx_lx_p}\phi\notag\\
	 	&+\xi\Delta x^2\partial_t\sum_{i,j=1}^d\nabla^2_{ij}\Big(\frac{1}{s_{x_i}}-1\Big)\Big(1-\frac{1}{s_{2|x_ix_j}}-\frac{1}{s_{x_i}}\Big)\bm{\varepsilon}_{2|x_ix_j}\phi+\sum_{i,j=1}^d\nabla^2_{ij}\Big(\frac{1}{s_{x_i}}-1\Big)\bm{\varepsilon}_{2|x_ix_j}\phi+\frac{\Delta x^2}{6}\sum_{i,j,l,p=1}^d\nabla^4_{ijlp}\Big(\frac{1}{s_{x_i}}-1\Big)\bm{\varepsilon}_{4|x_ix_jx_lx_p}\phi \notag\\
	 	&+\Delta x^2\sum_{i,j,l,p=1}^d\nabla^4_{ijlp}\Big(\frac{1}{s_{3|x_ix_jx_l}s_{2|x_ix_j}}-\frac{1}{2s_{3|x_ix_jx_l}}-\frac{1}{2s_{2|x_ix_j}}\Big)\Big(\frac{1}{s_{x_i}}-1\Big)\bm{\varepsilon}_{4|x_ix_jx_lx_p}\phi\notag\\
	 	&  
	 	+\xi\Delta x^2\sum_{i,j=1}^d\nabla^2_{ij}\Big(\frac{1}{s_{2|x_ix_j}}-\frac{1}{2}\Big)\Big(\frac{1}{s_{x_i}}-\frac{1}{2}\Big)\eta\bm{\varepsilon}_{2|x_ix_j}\phi	 \notag\\
	 	&+\frac{1}{2}\xi\Delta x^2\sum_{i,j=1}^d \nabla^2_{ij}\Big(1-\frac{1}{s_{2|x_ix_j}}\Big)\partial_t\bm{\varepsilon}_{2|x_ix_j}  -\frac{1}{2}\Delta x^2 \sum_{i,j,l,p=1}^d\nabla^4_{ijlp}\Big(1-\frac{1}{s_{2|x_ix_j}}\Big)\Big(\frac{1}{s_{3|x_ix_jx_l}}-\frac{1}{2}\Big)\bm{\varepsilon}_{4|x_ix_jx_lx_p}\phi \notag\\
	 	&-\frac{\Delta x^2}{6}\sum_{i,j,l,p=1}^d\nabla^4_{ijlp}\Big(1-\frac{1}{s_{3|x_ix_jx_l}}\Big)\bm{\varepsilon}_{4|x_ix_jx_lx_p}\phi 	+\xi R
.  
	 \end{align} 
	 	 From the above equation and  with respect to the diagonal-anisotropic diffusion equation (\ref{DE}), the zeroth-order  truncation error ($\makebox{Err}_{\makebox{0th}}$) of the MRT-LB model (\ref{lb}) is expressed as 
	 \begin{align} \label{err0}
	 	\makebox{Err}_{\makebox{0th}}=\Big(\frac{1}{s_{x_i}}-\frac{1}{2}\Big)\big(2\omega_i+4(d-1)\tilde{\omega}\big) +\frac{\Delta t}{s_{x_i}} \Big(1-\frac{1}{s_{x_i}}\Big)\eta\big(2\omega_i+4(d-1)\tilde{\omega}\big)-\kappa_{x_i}\xi,\quad i\in\llbracket 1,d\rrbracket, 
	 \end{align}   
	 where (\ref{moment}) and the fact of $\xi\partial_t\phi=\sum_{i,j=1}^d\nabla^2_{ij}(1/s_{x_i}-1/2)\bm{\varepsilon}_{2|x_ix_j}\phi+\xi R+O(\Delta x^2)$ have been used. Here, we would like to point out that if we introduce a notation $\tilde{s}_{x_i}$  satisfying
	 \begin{align*}
	 	\tilde{\epsilon}_{x_i}=\epsilon_{x_i}+\frac{\Delta t}{s_{x_i}} \Big(1-\frac{1}{s_{x_i}}\Big)\eta\big(2\omega_i+4(d-1)\tilde{\omega}\big),
	 \end{align*}
	 here    $\tilde{\epsilon}_{x_i}=(1/\tilde{s}_{x_i}-1/2)\big(2\omega_i+4(d-1)\tilde{\omega}\big)$ and $\epsilon_{x_i}=(1/s_{x_i}-1/2)\big(2\omega_i+4(d-1)\tilde{\omega}\big)$, one can obtain
	 \begin{align}\label{err0-sol}
	 	s_{x_i}=\frac{2\eta\Delta t}{\Delta t\eta-\sqrt{\big(\tilde{s}_{x_i}-4\Delta t\eta+\Delta t^2\tilde{s}_{x_i}\eta^2+2\Delta t\tilde{s}_{x_i}\eta\big)/\tilde{s}_{x_i}}+1}=\tilde{s}_{x_i}+O(\Delta x^2).
	 \end{align}
	 Moreover, if    (\ref{err0}) holds, i.e., $\kappa_{x_i}=\tilde{\epsilon}_i/\xi$,  we  have
	 \begin{subequations}\label{seoncd-order-relation}
	 	 \begin{align}
	 		&\partial_t\phi=\sum_{i=1}^d\kappa_{x_i}\nabla_i^2\phi+R+O(\Delta x^2),\\
	 		&\partial_t^2\phi=\sum_{i,j=1}^d\kappa_{x_i}\kappa_{x_j}\nabla_{iijj}^4\phi +2\eta\sum_{i=1}^d\kappa_{x_i}\nabla_i^2\phi+\eta R+O(\Delta x^2).
	 	\end{align}
	 \end{subequations}

	 In the sense of fourth-order accuracy, substituting (\ref{seoncd-order-relation}) and  (\ref{err0-sol}) into (\ref{err}) give rises to
	 \begin{align}\label{macro-simp}
	 	&\partial_t\phi+\frac{\xi\Delta x^2}{2}\sum_{i,j=1}^d\kappa_{x_i}\kappa_{x_j} \nabla_{iijj}^4\phi +\frac{\Delta x^2}{2}\sum_{i=1}^d\Big(\frac{1}{s_{x_i}}-\frac{1}{2}\Big)\bm{\varepsilon}_{2|x_i^2}\eta\nabla_i^2\phi +O(\Delta x^4) =  \sum_{i,j=1}^d\nabla^2_{ij}\frac{1}{\xi}\Big(\frac{1}{s_{x_i}}-\frac{1}{2}\Big)\bm{\varepsilon}_{2|x_ix_j}\phi\notag\\
	 	&
	 	+\frac{\Delta x^2}{\xi}\sum_{i,j,l,p=1}^d\Bigg(-\frac{7}{24}+\frac{1}{6s_{x_i}}+\Big(\frac{1}{\tilde{s}_{3|x_ix_jx_l}s_{2|x_ix_j}}-\frac{1}{2\tilde{s}_{3|x_ix_jx_l}}-\frac{1}{2s_{2|x_ix_j}}\Big)\Big(\frac{1}{s_{x_i}}-1\Big)\notag\\
	 	&\qquad\qquad\qquad\qquad
	 	-\Big(\frac{1}{2}-\frac{1}{2s_{2|x_ix_j}}\Big)\Big(\frac{1}{\tilde{s}_{3|x_ix_jx_l}}-\frac{1}{2}\Big)+\frac{1}{6\tilde{s}_{3|x_ix_jx_l}}\Bigg)\nabla^4_{ijlp}\bm{\varepsilon}_{4|x_ix_jx_lx_p}\phi\notag\\
	 	&+\Delta x^2 \sum_{i,j=1}^d\nabla^4_{iijj}\bigg(\Big(\frac{1}{\tilde{s}_{x_i}}-1\Big)\Big(1-\frac{1}{s_{2|x_i^2}}-\frac{1}{\tilde{s}_{x_i}}\Big)+\frac{1}{2}-\frac{1}{2s_{2|x_i^2}}\bigg)\kappa_{x_j}\bm{\varepsilon}_{2|x_i^2}\phi  \notag\\ 
	 	&+\Delta x^2 \sum_{i=1}^d\nabla^2_{ii}\bigg(\Big(\frac{1}{s_{x_i}}-1\Big)\Big(1-\frac{1}{s_{2|x_i^2}}-\frac{1}{s_{x_i}}\Big)+\frac{1}{2}-\frac{1}{2s_{2|x_i^2}}\bigg)\eta\bm{\varepsilon}_{2|x_i^2}\phi   +R 	+\Delta x^2\sum_{i=1}^d\nabla^2_{ii}\Big(\frac{1}{s_{2|x_i^2}}-\frac{1}{2}\Big)\Big(\frac{1}{s_{x_i}}-\frac{1}{2}\Big)\eta\bm{\varepsilon}_{2|x_i^2}\phi, 
	 \end{align} 
	  where \begin{align}\label{s3}
	   \tilde{s}_{3|x_kx_lx_p}= \begin{cases}
	  		\tilde{s}_{x_k},& k\in\llbracket 1, d\rrbracket,\quad p=l=k,\\
	  		\tilde{s}_{3|x_k^2x_l},& k,l\in\llbracket 1,d\rrbracket,\quad p=k,\:\: k<l.
	  	\end{cases} 
	  \end{align} 
	
From Eq. (\ref{macro-simp}), one can obtain the following second-order  truncation error ($\makebox{Err}_{\makebox{2rd}}$) of the MRT-LB model,
	\begin{align} \label{err2}
		\makebox{Err}_{\makebox{2rd}}= \Delta x^2\Bigg(& -\frac{\xi}{2}\sum_{i,j=1}^d\kappa_{x_i}\kappa_{x_j}\nabla_{iijj}^4 \phi+ \frac{1}{\xi}\sum_{i=1 }^d\tilde{\makebox{err}}^{(1)}_{x_i|x_i^2|x_i^3}\nabla^4_{i}\big(2\omega_i+4(d-1)\tilde{\omega}\big)\phi\notag\\
		&+\frac{4}{\xi}\sum_{i,j=1(i\neq j)}^d\Big(\tilde{\makebox{err}}^{(1)}_{x_i|x_i^2|x_i^2x_j}+\tilde{\makebox{err}}^{(1)}_{x_i|x_ix_j|x_i^2x_j}+\tilde{\makebox{err}}^{(1)}_{x_i|x_ix_j|x_ix_j^2}\Big)\nabla^4_{iijj} \tilde{\omega}\phi\notag\\
		&+ \bigg(  \sum_{i,j=1(i\neq j)}^d\nabla^4_{iijj} \tilde{\makebox{err}}^{(2)}_{{x_i}}\kappa_{x_j}+\sum_{i=1}^d\nabla^4_{i}\tilde{\makebox{err}}^{(2)}_{{x_i}}\kappa_{x_i} +  \sum_{i=1}^d\nabla^2_{i}\frac{\eta}{s_{x_i}} \Big(1-\frac{1}{s_{x_i}}\Big)\bigg)\big(2\omega_i+4(d-1)\tilde{\omega}\big)\phi \Bigg),  
	\end{align}   
where (\ref{moment}) and (\ref{la-m-s})
 have been used, and 
 \begin{align*}
	&\tilde{\makebox{err}}^{(1)}_{{\alpha}|{\beta}|{\gamma}}=-\frac{7}{24}+\frac{1}{6\tilde{s}_{{\alpha}}}+\Big(\frac{1}{s_{3|\gamma}s_{2|\beta}}-\frac{1}{2s_{3|\gamma}}-\frac{1}{2s_{2|\beta}}\Big)\Big(\frac{1}{\tilde{s}_{\alpha}}-1\Big)
	-\Big(\frac{1}{2}-\frac{1}{2s_{2|\beta}}\Big)\Big(\frac{1}{s_{3|\gamma}}-\frac{1}{2}\Big)+\frac{1}{6s_{3|\gamma}},\\
	&\tilde{\makebox{err}}^{(2)}_{{\alpha}}=\Big(\frac{1}{\tilde{s}_{\alpha}}-1\Big)\Big(1-\frac{1}{s_{2|\alpha^2}}-\frac{1}{\tilde{s}_{\alpha}}\Big)+\frac{1}{2}-\frac{1}{2s_{2|\alpha^2}}.
\end{align*}

	Thus, according to (\ref{err0}), (\ref{err0-sol}), and (\ref{err2}), in order to ensure that the modified equation (\ref{macro-simp}) of the MRT-LB model (\ref{lb}) can be fourth-order consistent with (\ref{DE}), the following conditions should be satisfied, 
	\begin{subequations}\label{fourth-conition}
		\begin{align}   
			&s_{x_i}=\frac{2\eta\Delta t}{\Delta t\eta-\sqrt{\big(\tilde{s}_{x_i}-4\Delta t\eta+\Delta t^2\tilde{s}_{x_i}\eta^2+2\Delta t\tilde{s}_{x_i}\eta\big)/\tilde{s}_{x_i}}+1},\quad i\in\llbracket 1,d\rrbracket,\label{fourth-coition-1}\\ 
			&\tilde{\epsilon}_{x_i}=\kappa_{x_i}\xi,\quad i\in\llbracket 1,d\rrbracket,\label{fourth-coition-2}\\ 
			&\frac{\tilde{\epsilon}_{x_i}}{2}\Big(\frac{1}{\tilde{s}_{x_i}}-\frac{1}{2}\Big) =\tilde{\makebox{err}}^{(1)}_{x_i|x_i^2|x_i^3}+  {\makebox{err}}^{(2)}_{x_i}\tilde{\epsilon}_{x_i},  \quad i\in\llbracket 1,d\rrbracket,\label{fourth-coition-3}\\ 
			&\tilde{\epsilon}_{x_i} \tilde{\epsilon}_{x_j}  =  4\Big(\tilde{\makebox{err}}^{(1)}_{x_i|x_i^2|x_i^2x_j}+\tilde{\makebox{err}}^{(1)}_{x_j|x_j^2|x_ix_j^2}
			+\tilde{\makebox{err}}^{(1)}_{x_i|x_ix_j|x_i^2x_j}+\tilde{\makebox{err}}^{(1)}_{x_j|x_ix_j|x_ix_j^2}
			+\tilde{\makebox{err}}^{(1)}_{x_i|x_ix_j|x_ix_j^2}+\tilde{\makebox{err}}^{(1)}_{x_j|x_ix_j|x_i^2x_j}\Big)\tilde{\omega} \notag\\
			&\qquad\qquad+ \tilde{\makebox{err}}^{(2)}_{x_i}\tilde{\epsilon}_{x_j}\big(2\omega_i+4(d-1)\tilde{\omega}\big) + \tilde{\makebox{err}}^{(2)}_{x_j}\tilde{\epsilon}_{x_i}\big(2\omega_j+4(d-1)\tilde{\omega}\big) ,\quad i,j\in\llbracket 1,d\rrbracket,\:\: i< j. \label{fourth-coition-4}\
		\end{align} 
	\end{subequations}
	And how to solve the  above $(d^2+5d)/2$ conditions can be found in the following Part \ref{automatic}. 
	
Similar to the above accuracy analysis, we can also derive the following fourth-order consistent conditions for the MRT-LB model (\ref{lb}) with the D$d$Q($2d+1$) lattice structure proposed in Remark 2, 
 \begin{subequations}\label{fourth-coition-ddq2d+1}
 	\begin{align}   
 		&s_{x_i}=\frac{2\eta\Delta t}{\Delta t\eta-\sqrt{\big(\tilde{s}_{x_i}-4\Delta t\eta+\Delta t^2\tilde{s}_{x_i}\eta^2+2\Delta t\tilde{s}_{x_i}\eta\big)/\tilde{s}_{x_i}}+1},\quad i\in\llbracket 1,d\rrbracket,\\
 		&\tilde{\epsilon}_{x_i}=\kappa_{x_i}\xi,\quad i\in\llbracket 1,d\rrbracket,\\ 
 		&\frac{\tilde{\epsilon}_{x_i}}{2}\Big(\frac{1}{\tilde{s}_{x_i}}-\frac{1}{2}\Big) = \frac{1}{3\tilde{s}_{x_i}}-\frac{7}{24}+\Big(\frac{1}{\tilde{s}_{x_i}s_{2|x_i^2}}-\frac{1}{2\tilde{s}_{x_i}}-\frac{1}{2s_{2|x_i^2}}\Big)\Big(\frac{1}{\tilde{s}_{x_i}}-1\Big)-\frac{1}{2}\Big(1-\frac{1}{s_{2|x_i^2}}\Big)\Big(\frac{1}{ \tilde{s}_{x_i}}-\frac{1}{2}\Big)+ \tilde{{\makebox{err}}}^{(2)}_{x_i}\tilde{\epsilon}_{x_i}  ,\quad i\in\llbracket 1,d\rrbracket, \\
 		&\tilde{\epsilon}_{x_i} \tilde{\epsilon}_{x_j}  =  \tilde{{\makebox{err}}}^{(2)}_{x_i}\tilde{\epsilon}_{x_i}\tilde{\epsilon}_{x_j}\big(2\omega_i \big) + \tilde{{\makebox{err}}}^{(2)}_{x_i}\tilde{\epsilon}_{x_j}\tilde{\epsilon}_{x_i}\big(2\omega_j \big),\quad i,j\in\llbracket 1,d\rrbracket,\:\: i< j,
 	\end{align} 
 \end{subequations}
 	 It can be easily verified that   (\ref{fourth-coition-ddq2d+1}) has no solution unless the diagonal-anisotropic diffusion equation (\ref{DE}) is reduced to the isotropic type (see the following Part \ref{automatic} for details). 
	\subsection{The fourth-order initialization schemes}\label{ini-sec}
	In order to obtain an overall fourth-order MRT-LB model for the diagonal-anisotropic diffusion equation (\ref{DE}), in this part, we discuss how to construct the fourth-order  initialization scheme.

		In fact, (\ref{4th}) can be used to initialize the unknown distribution function \cite{Chen2023}. While in this work, we will adopt the following scheme to approximate the distribution function at the initial time, i.e.,
		\begin{align}\label{ini}
			\mathbf{f}(\mathbf{x},0)= \mathbf{f}^{eq}(\mathbf{x},0)-\Delta x\bm{\Lambda}^{-1}\mathbf{e}_i\cdot\nabla\mathbf{f}^{eq}(\mathbf{x},0),
		\end{align}
		where $f^{eq}_k(\mathbf{x},0)=\omega_k\phi^0(\mathbf{x})$. Although (\ref{2rd}) shows that (\ref{ini}) is only a second-order approximation of the distribution function $f_k$, we will prove below that it is actually a fourth-order  method.  We first substitute (\ref{ini})  into the MRT-LB model (\ref{lb}), and obtain
		\begin{align}\label{ini-dt}
			\mathbf{f}(\mathbf{x},\Delta t)=\mathbf{T}\bigg(\mathbf{f}^{eq}-\Delta x\bm{\Lambda}^{-1}\mathbf{e}_i\cdot\nabla\mathbf{f}^{eq} +\Delta t\Big(\mathbf{I}-\frac{\bm{\Lambda}}{2}\Big)\mathbf{R}\bigg)(\mathbf{x},0), 
		\end{align}
		where $\mathbf{T}$ is the space shift matrix  defined by
		\begin{align*}
			\mathbf{T}=\makebox{\textbf{diag}}\big(T^{-\mathbf{e}_1}_{\Delta x},T^{-\mathbf{e}_2}_{\Delta x},\cdots,T^{-\mathbf{e}_q}_{\Delta x}\big) \makebox{ with } T^{-\mathbf{e}_i}_{\Delta x}f_k(\mathbf{x},t)=f_k(\mathbf{x}-\mathbf{e}_i\Delta x,t),\quad i\in\llbracket 1,q\rrbracket.
		\end{align*}
		After taking the zeroth moment of $f_k$ in (\ref{ini-dt}), we have
		\begin{align*} 
			\sum_{k=1}^qf_k(\mathbf{x},\Delta t)=&\phi(\mathbf{x},0)+\Delta t\frac{\Delta x^2}{\Delta t}\sum_{i=1}^d\big(2\omega_i+4(d-1)\tilde{\omega}\big)\Big(\frac{1}{s_{x_i}}-\frac{1}{2}\Big)\nabla^2_{i}\phi(\mathbf{x},0)+\frac{\Delta t}{2}R(\mathbf{x},0)+O(\Delta x^2)\notag\\
			=& \phi(\mathbf{x},0)-\frac{\Delta t}{2}R(\mathbf{x},0)+\Delta t\partial_t\phi +O(\Delta x^4)=\phi(\mathbf{x},\Delta t)-\frac{\Delta t}{2}R(\mathbf{x},\Delta t)+O(\Delta x^4), 
		\end{align*}
		where the fact of $T^{-\mathbf{e}_i}_{\Delta x}=\exp(                                                                                                                                                                                                                              -\mathbf{e}_i\cdot\nabla \Delta x)$ has been used. Thus, the approximation of $\phi(\mathbf{x},\Delta t)$ through adopting the initialization scheme (\ref{ini}) can achieve a fourth-order accuracy. 
		
		We now further focus on the time $t=2\Delta t$. Substituting (\ref{ini-dt}) into (\ref{lb}) gives rise to
		\begin{align}\label{ini-dt2}
			\mathbf{f}(\mathbf{x},2\Delta t)=&\mathbf{T}\bigg(\mathbf{f}-\bm{\Lambda}(\mathbf{f}-\mathbf{f}^{eq})+\Delta t\Big(\mathbf{I}-\frac{\bm{\Lambda}}{2}\Big)\mathbf{R}\bigg)(\mathbf{x},\Delta t)\notag\\
			=&\mathbf{T}\bigg(\mathbf{I}-\bm{\Lambda}(\mathbf{I}-\mathbf{E})+\Delta t\Big(\mathbf{I}-\frac{\bm{\Lambda}}{2}\Big)\eta\mathbf{E}\bigg)\mathbf{T}\bigg(\mathbf{f}^{eq}-\Delta x\bm{\Lambda}^{-1}\mathbf{e}_i\cdot\nabla\mathbf{f}^{eq} +\Delta t\Big(\mathbf{I}-\frac{\bm{\Lambda}}{2}\Big)\mathbf{R}\bigg)(\mathbf{x},0)+\Delta t\Big(\mathbf{I}-\frac{\bm{\Lambda}}{2}\Big)\mathbf{S}_R\notag\\
			=&\Delta t\Big(\mathbf{I}-\frac{\bm{\Lambda}}{2}\Big)\mathbf{S}_R+\Big(\mathbf{I}-\Delta x\mathbf{D}+\frac{\Delta x^2\mathbf{D}}{2}-\frac{\Delta x^3\mathbf{D}^3}{6}+O(\Delta x^4)\Big)\bigg(\mathbf{I}-\bm{\Lambda}(\mathbf{I}-\mathbf{E})+\Delta t\Big(\mathbf{I}-\frac{\bm{\Lambda}}{2}\Big)\eta\mathbf{E}\bigg)\notag\\
			&\Big(\mathbf{I}-\Delta x\mathbf{D}+\frac{\Delta x^2\mathbf{D}}{2}-\frac{\Delta x^3\mathbf{D}^3}{6}+O(\Delta x^4)\Big)\bigg(\mathbf{f}^{eq}-\Delta x\bm{\Lambda}^{-1}\mathbf{e}_i\cdot\nabla\mathbf{f}^{eq} +\Delta t\Big(\mathbf{I}-\frac{\bm{\Lambda}}{2}\Big)\mathbf{R}\bigg)(\mathbf{x},0),
		\end{align}
		where  $\mathbf{D}$, $\mathbf{E}$, and $\mathbf{S}_R$ are given by
		\begin{align*}
			&\mathbf{D}=\makebox{\textbf{diag}}\big(\mathbf{e}_1\cdot \nabla, \mathbf{e}_2\cdot \nabla,\cdots,\mathbf{e}_q\cdot \nabla\big),\\& \mathbf{E}=\big(\omega_0\mathbf{I}_{q\times 1},\omega_1\mathbf{I}_{q\times 1},\cdots \omega_{q-1}\mathbf{I}_{q\times 1}\big)^T,\\& \mathbf{S}_R=\mathbf{M}^{-1}\Big(1,\mathbf{0}_{d\times 1}, \bm{\varepsilon}_{2|\mathbf{x}^{.2}},\mathbf{0}_{d(d-1)/2},\mathbf{0}_{d(d-1)}, \bm{\varepsilon}_{4|(\mathbf{x}^{.2}\mathbf{x}^{.2})_{\alpha<\beta}}\Big)^TS.
		\end{align*}
		After some algebraic manipulations, the zeroth-order moment of the above equation (\ref{ini-dt2}) is expressed as
		\begin{align}\label{ini-2dt}
			\sum_{k=1}^qf_k(\mathbf{x},2\Delta t)=\phi(\mathbf{x},0)-\frac{\Delta t}{2}R(\mathbf{x},0)+2\Delta t\partial_t\phi(\mathbf{x},0)+O(\Delta x^4)=\phi(\mathbf{x},2\Delta t)-\frac{\Delta t}{2}R(\mathbf{x},2\Delta t)+O(\Delta x^4),
		\end{align}
		where the fourth-order consistent conditions in (\ref{fourth-conition})  have been used, which  indicates that the variable $\phi$ at the time $t=2\Delta t$ can be calculated by $\sum_{k=1}^qf_k-\Delta tR/2$ with a fourth-order accuracy. And the analysis of the time $t=n\Delta t$ ($n> 2$) is similar, we do not present the details here.
	 
	\subsection{The stability analysis of the MRT-LB model}\label{stability}
	In this part, we will discuss whether the MRT-LB model (\ref{lb})  can satisfy the stability structure, and further analyze the $L^2$ stability.  Without loss of generality, in the following analysis we neglect the constant source term of the diagonal-anisotropic diffusion equation (\ref{DE}).  
	
	Firstly, we present the following definition of the stability structure \cite{Banda2006,JunkYang2009}.\\
	\textbf{Definition 1.} Set $\mathbf{Q}(\mathbf{f})=\bm{\Lambda}(\mathbf{f}^{eq}-\mathbf{f})+\Delta t\big(\mathbf{I}-\bm{\Lambda}/2\big)\mathbf{R}$ and $\mathbf{J}(\mathbf{f})=\mathbf{Q}'(\mathbf{f})$. Then let $\mathbf{f}=\tilde{\mathbf{f}}$ which satisfies $\mathbf{Q}(\tilde{\mathbf{f}})=\mathbf{0}$. If there exists an invertible matrix $\mathbf{P}\in\mathbb{R}^{q\times q}$ such that $\mathbf{P}^T\mathbf{P}=\makebox{\textbf{diag}}(a_1,a_2,\cdots,a_q)$ ($a_k> 0$) is diagonal and 
	\begin{align}
		\mathbf{P}\mathbf{J}(\tilde{\mathbf{f}})=-\makebox{\textbf{diag}}(\lambda_1,\lambda_2,\cdots,\lambda_q)\mathbf{P},
	\end{align}
	where $\lambda_k>0$ for all $k$, we claim  that the LB method satisfies the stability structure $(\mathbf{P},\mathbf{a},\bm{\lambda})$ at $\mathbf{f}=\tilde{\mathbf{f}}$.
	
Based on the previous work \cite{Yang2024}, which provides an automatic approach to determine the condition that the LB method can satisfy the stability structure, one can decompose $\mathbf{J}$ as
	\begin{align}\label{J}
		\mathbf{J}=\big(\mathbf{J}\mathbf{W}\big)\mathbf{W}^{-1} \makebox{  with  } \mathbf{W}=\makebox{\textbf{diag}}(\omega_0,\omega_1,\cdots,\omega_{q-1}).
	\end{align}
Evidently, if matrix $\big(\mathbf{J}\mathbf{W}\big)$ is symmetric and semi-negative definite,  $\big(\mathbf{W}^{-1/2}\mathbf{JW}\mathbf{W}^{-1/2}\big)$ still remains the symmetry and semi-negative definiteness. Thus, there must exit an orthogonal  matrix $\mathbf{U}$ and a diagonal matrix $\makebox{\textbf{diag}}(\lambda_1,\lambda_2,\cdots,\lambda_q)$ such that 
	\begin{align*}
		\mathbf{W}^{-1/2}\big(\mathbf{J}\mathbf{W}\big)\mathbf{W}^{-1/2}=-\mathbf{U}^T\makebox{\textbf{diag}}(\lambda_1,\lambda_2,\cdots,\lambda_q)\mathbf{U},
	\end{align*}
then, (\ref{J}) can be rewritten as
\begin{align*}
	\mathbf{J}=\big(\mathbf{J}\mathbf{W}\big)\mathbf{W}^{-1} =-\mathbf{W}^{1/2}\mathbf{U}^T\makebox{\textbf{diag}}(\lambda_1,\lambda_2,\cdots,\lambda_q)\mathbf{U}\mathbf{W}^{1/2}\mathbf{W}^{-1}.
\end{align*}
In particular, let $\mathbf{P}=\mathbf{U}\mathbf{W}^{-1/2}$, we have
\begin{align*}
	\mathbf{P}^T\mathbf{P}=\mathbf{W}^{-1/2}\mathbf{U}^T\mathbf{U}\mathbf{W}^{-1/2}=\mathbf{W}^{-1},
\end{align*}
which is consistent with Definition 1. This means that to analyze whether the LB model  satisfies the stability structure, we only need to find the condition that ensures the matrix $(\mathbf{JW})$ to be  symmetric and semi-negative definite. And it is noteworthy that if we further set $\lambda_k$, which corresponds to the non-conservative moment, to be in the range of $(0,2)$, the convergence of the LB method  can also be proved \cite{JunkYong2009,junk2008convergence,JunkYang2009}.

For the present MRT-LB model (\ref{lb}), the matrix $\mathbf{J}$ is given by
\begin{align}
	\mathbf{J}=\bm{\Lambda}\big(\mathbf{E}-\mathbf{I}\big)+\frac{2\Delta t\eta}{2-\Delta t\eta}\mathbf{E}.
\end{align}
To determine the condition that ensures the matrix $(\mathbf{JW})$ to be symmetric and semi-negative definite, the derivation process is divided into the subsequent two steps.

\begin{itemize}
	\item Symmetry
	 
	It is obvious that  $(\mathbf{EW})$ is a symmetric matrix, thus we just need to find the condition that ensures the matrix $\big(\bm{\Lambda}(\mathbf{E}-\mathbf{I})\mathbf{W}\big)$ to be symmetric, which can be derived by a simple computer code and is given by
	\begin{align}\label{sym}
		s_{3|x_i^2x_j}=s_{x_j},\quad  s_{4|x_i^2x_j^2}=s_{2|x_i^2}=s_{2|x_j^2},\quad i,j\in\llbracket1,d\rrbracket,\:\: i<j.
	\end{align}
	\item Semi-negative definiteness
	
	The eigenvalues $\lambda_k$ ($k\in\llbracket 1,q\rrbracket$) of matrix $\mathbf{J}$ are $2\Delta t\eta/(2-\Delta \eta)$ (corresponds to the conservative moment $\phi$), $-s_{x_i}$, $-s_{2|x_i^2}$, $-s_{2|x_ix_j}$, $-s_{3|x_i^2x_j}$, and $-s_{4|x_i^2x_j^2}$. Therefore, provided that the condition (\ref{condition-stability}) is satisfied and all the relaxation parameters, except $s_0$, are located in the range of $(0,2)$, i.e.,
	\begin{align}\label{condition-stability}
		s_{2|x_ix_j}\in(0,2),\quad s_{3|x_i^2x_j}=s_{x_j}\in(0,2),\quad  s_{4|x_i^2x_j^2}=s_{2|x_i^2}=s_{2|x_j^2}\in(0,2),\quad i,j\in\llbracket1,d\rrbracket,\:\: i<j,
	\end{align}
	the matrix $(\mathbf{JW})$ is semi-negative definiteness. What is more, the convergence of the present MRT-LB model (\ref{lb}) can be proved \cite{JunkYong2009,junk2008convergence,JunkYang2009}.
\end{itemize}

\begin{figure} [H] 
	\begin{center}  
		\subfloat[$\tilde{\omega}=1/18$]
		{
			\includegraphics[width=0.3\textwidth]{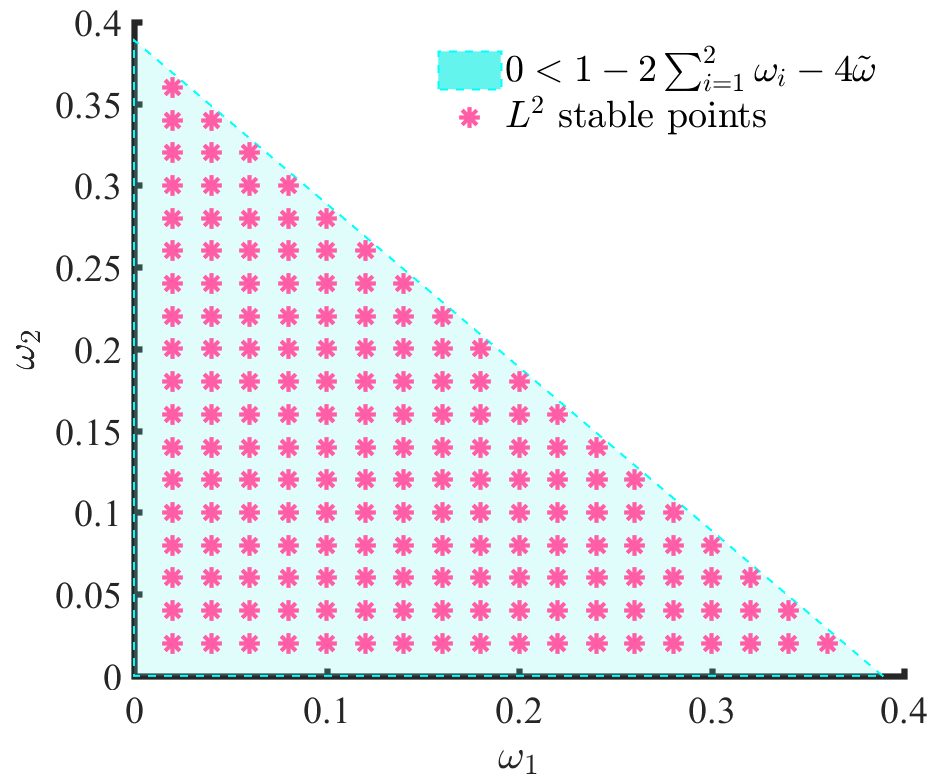}
		} 
		\subfloat[$\tilde{\omega}=1/36$]
		{
			\includegraphics[width=0.3\textwidth]{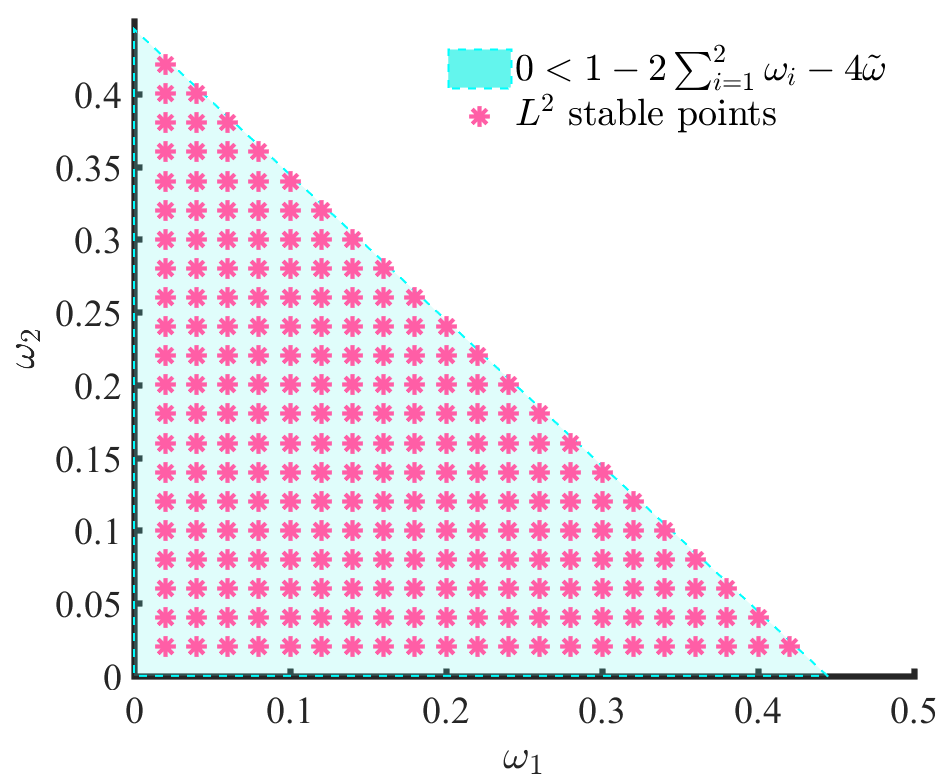}
		}  
		\subfloat[$\tilde{\omega}=1/180$]
		{
			\includegraphics[width=0.3\textwidth]{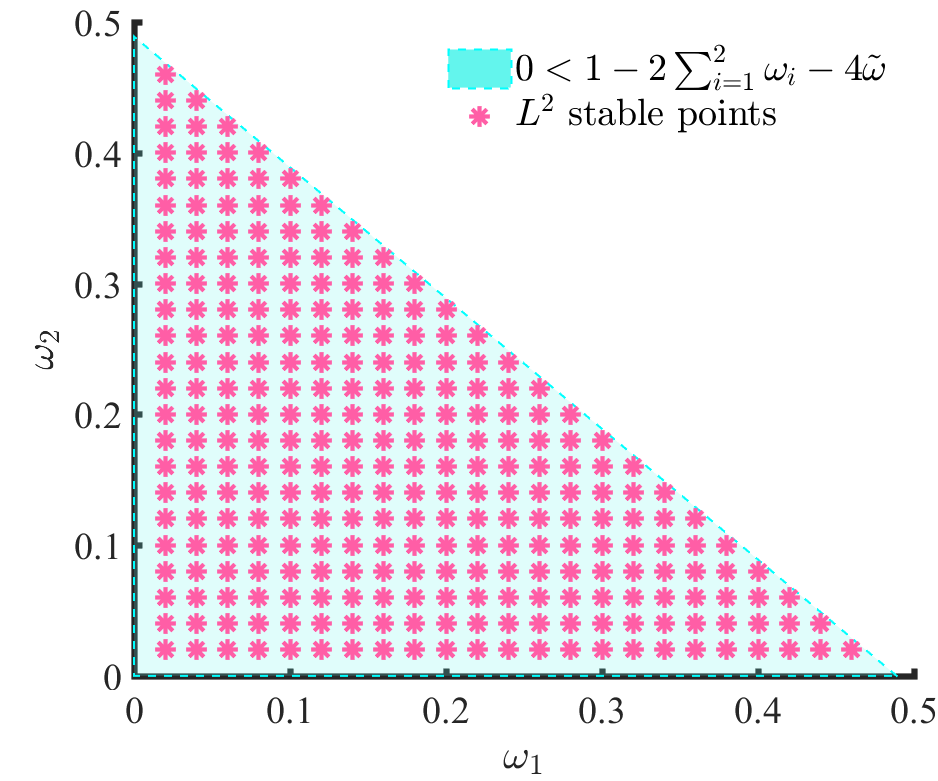}
		}  
		
		\subfloat[$\tilde{\omega}=1/18$]
		{
			\includegraphics[width=0.3\textwidth]{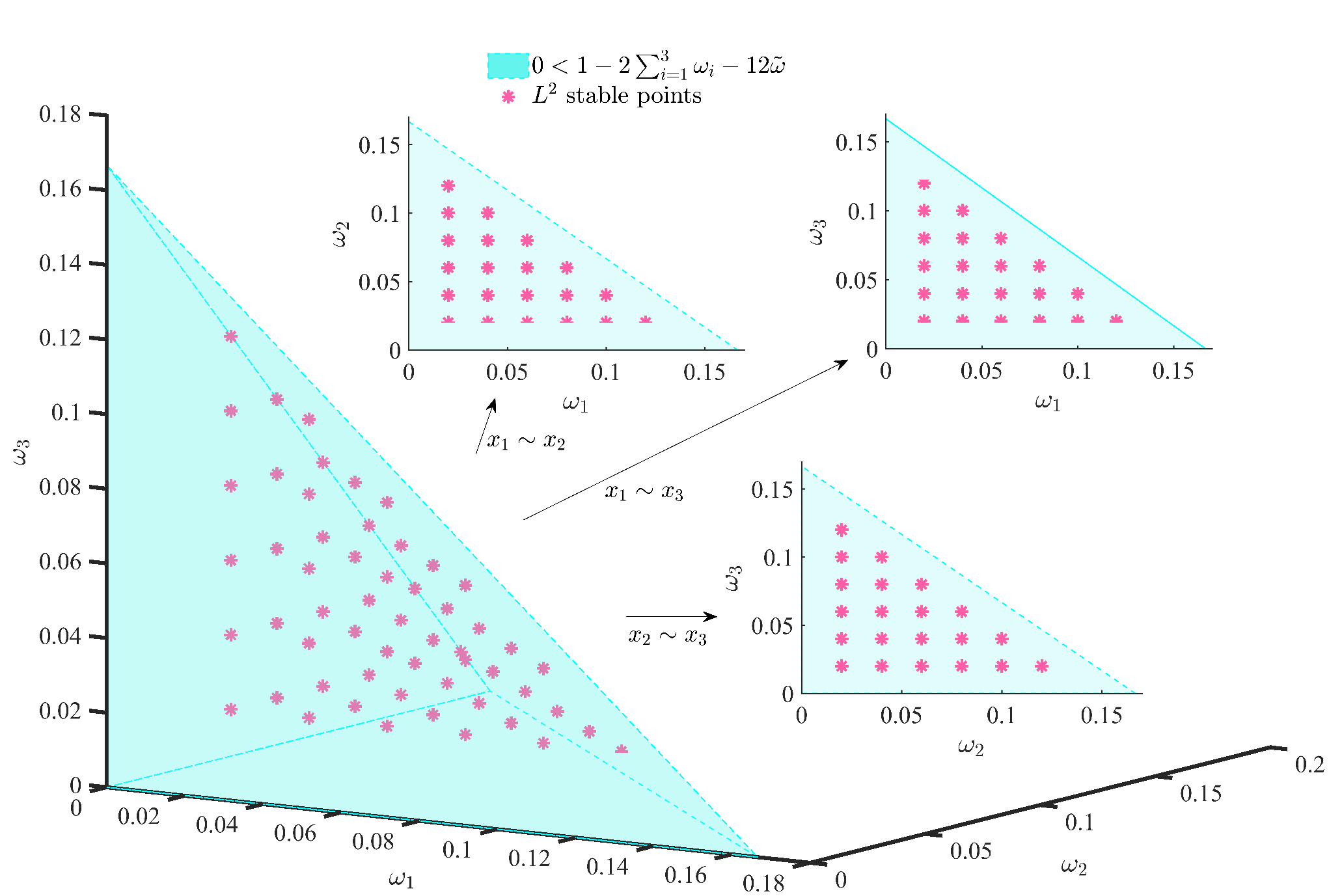}
		} 
		\subfloat[$\tilde{\omega}=1/36$]
		{
			\includegraphics[width=0.3\textwidth]{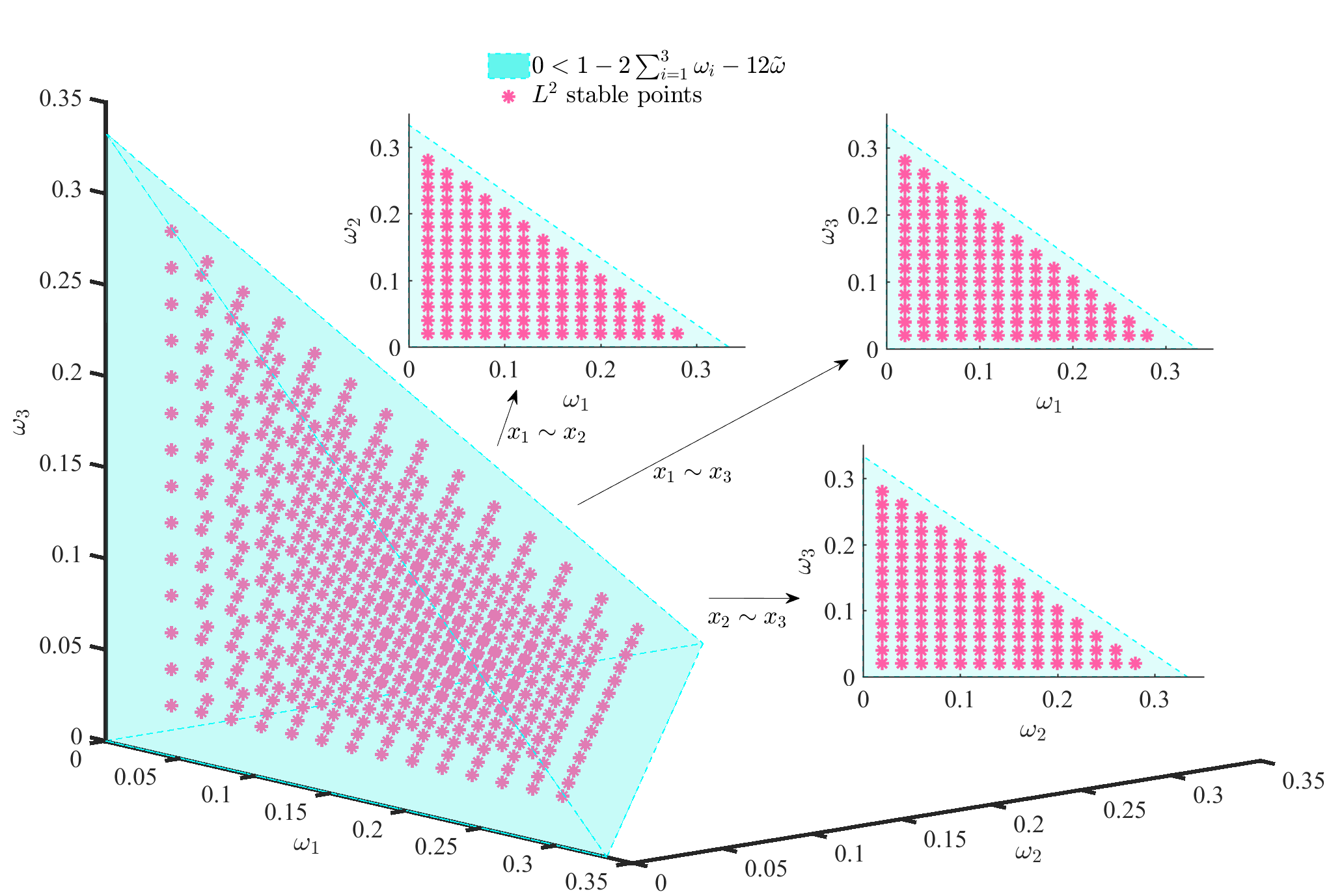}
		}  
		\subfloat[$\tilde{\omega}=1/180$]
		{
			\includegraphics[width=0.3\textwidth]{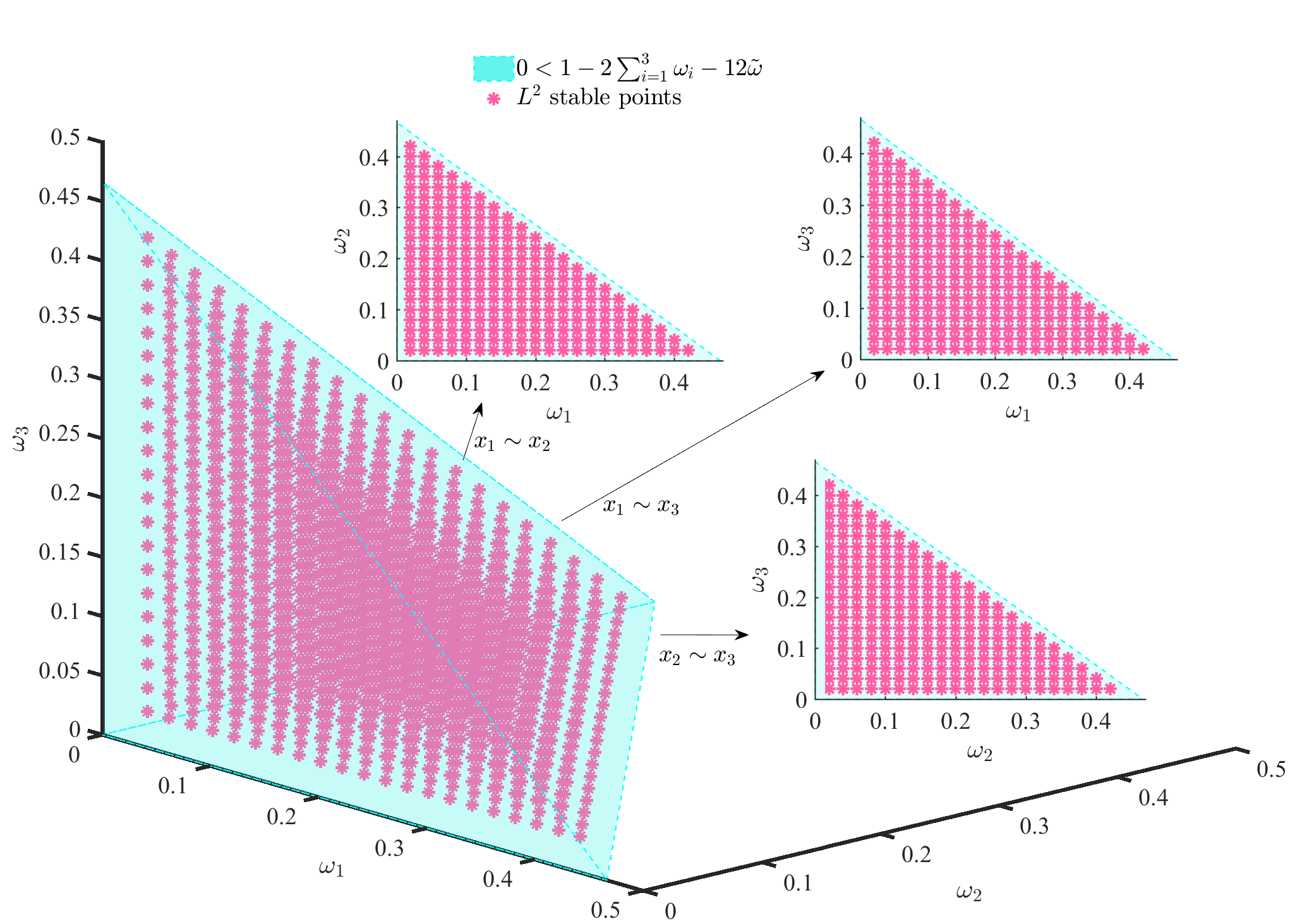}
		}  
		\caption{$L^2$ stability regions of the MRT-LB model for the two-dimensional (top) and three-dimensional (bottom) cases under  different values of the weight coefficient $\tilde{\omega}$.} \label{fig-region-stability}\end{center}  
\end{figure}

Based on the analysis results of the structure stability above, we now conduct an analysis on the $L^2$ stability of the MRT-LB model (\ref{lb}). To begin with, we introduce a weight $L^2$-norm \cite{JunkYang2009}
\begin{align*}
	\|\mathbf{f}(\mathbf{x}_i,t_n)\|_{\mathbf{P}}=\Big|\mathbf{Pf}(\mathbf{x}_i,t_n)\Big|=\bigg(\sum_{k=1}^qa_kf_k^2(\mathbf{x}_i,t_n)\bigg)^{1/2},
\end{align*}
and the associated operator norm is then given by
\begin{align*}
	\|\mathbf{B}\|_{\mathbf{P}}=\sup_{\mathbf{f}\neq \mathbf{0}}\frac{\|\mathbf{Bf}\|_{\mathbf{P}}}{\|\mathbf{f}\|_{\mathbf{P}}},
\end{align*}
Similarly, for matrix $\mathbf{C}\in\mathbb{C}^{q\times q}$, we define
\begin{align}\label{matrix-norm}
	\|\mathbf{C}\|_{\mathbf{P}}=\sup_{\mathbf{f}\neq \mathbf{0}}\frac{\|\mathbf{Cf}\|_{\mathbf{P}}}{\|\mathbf{f}\|_{\mathbf{P}}}.
\end{align}
For the present MRT-LB model (\ref{lb}), the amplification matrix $\mathbf{G}$ is given by
\begin{align}
	\mathbf{G}=\mathbf{\hat{T}}(\mathbf{I}+\mathbf{J})\makebox{ with }\mathbf{\hat{T}}=\makebox{\textbf{diag}}\Big(e^{-\mathbf{e}_1\cdot\bm{\xi}},e^{-\mathbf{e}_2\cdot\bm{\xi}},\cdots,e^{-\mathbf{e}_q\cdot\bm{\xi}}\Big),\quad \bm{\xi}\in [0,2\pi)^d,
\end{align}
then according to the definition of the matrix norm in (\ref{matrix-norm}), we can obtain 
\begin{align}\label{g-neq1}
	\|\mathbf{G}\|_{\mathbf{P}}=&\sup_{\mathbf{f}\neq \mathbf{0}}\frac{\|\mathbf{Gf}\|_{\mathbf{P}}}{\|\mathbf{f}\|_{\mathbf{P}}}=\sup_{\mathbf{f}\neq \mathbf{0}}\frac{\Big(\mathbf{f}^T\big(\mathbf{I}+\mathbf{J}\big)^T\mathbf{\hat{T}}^{*}\mathbf{P}^T\mathbf{P\hat{T}}\big(\mathbf{I}+\mathbf{J}\big)\mathbf{f}\Big)^{1/2}}{\Big(\mathbf{f}^T\mathbf{P}^T\mathbf{P}\mathbf{f}\Big)^{1/2}}=\sup_{\mathbf{f}\neq \mathbf{0}}\frac{\Big(\mathbf{f}^T\big(\mathbf{I}+\mathbf{J}\big)^T\mathbf{P}^T\mathbf{\hat{T}}^{*}\mathbf{\hat{T}P}\big(\mathbf{I}+\mathbf{J}\big)\mathbf{f}\Big)^{1/2}}{\Big(\mathbf{f}^T\mathbf{P}^T\mathbf{P}\mathbf{f}\Big)^{1/2}}\notag\\
	=&\sup_{\mathbf{f}\neq \mathbf{0}}\frac{\Big(\mathbf{Pf}^T\makebox{\textbf{diag}}\big((1-\lambda_1)^2,(1-\lambda_2)^2,\cdots,(1-\lambda_q)^2\big)\mathbf{Pf}\Big)^{1/2}}{\Big(\mathbf{f}^T\mathbf{P}^T\mathbf{P}\mathbf{f}\Big)^{1/2}}\notag\\
	=&\sup_{\mathbf{f}\neq \mathbf{0}}\frac{\bigg(\sum_{i=1}^q(1-\lambda_i)^2\Big(\sum_{j=1}^q\mathbf{P}_{ij}f_j\Big)\bigg)^{1/2}}{\bigg(\sum_{i=1}^q\Big(\sum_{j=1}^q\mathbf{P}_{ij}f_j\Big)\bigg)^{1/2}}=1,
\end{align} 
thus, from (\ref{g-neq1}), the spectral radius of matrix $\mathbf{G}$ satisfies
\begin{align}\label{rho-G}
	\rho(\mathbf{G})\leq \|\mathbf{G}\|_{\mathbf{P}}=1.
\end{align}
From (\ref{rho-G}), it is obvious that the characteristic polynomial of the matrix $\mathbf{G}$ [hereafter we denote it as $p(\mathbf{G})$] is a von Neumann polynomial. However, it remains unknown whether $p(\mathbf{G})$ is a simple von Neumann polynomial, which is the necessary and sufficient condition of the $L^2$ stability of the LB method  \cite{Chen2022,Chen2023}, yet it is generally challenging  to analyze theoretically. In this work, we will examine whether the roots lying on the unit circle are simple from a numerical  point of view.  For this purpose, we first set the relaxation parameters to satisfy the stability structure preserving condition (\ref{condition-stability}), then we focus on the $L^2$ stability region with respect to the weight coefficient $\omega_i$. The results for the most commonly used two- and three-dimensional cases are presented in Fig. \ref{fig-region-stability}, where different values of the weight coefficient $\tilde{\omega}$ are considered.  From Fig. \ref{fig-region-stability},  it can be observed that $p(\mathbf{G})$ is indeed a simple von Neumann polynomial unless the weight coefficient $\omega_i$ is not located in the range of (0,1), which is why the $L^2$ stability region in Fig. \ref{fig-region-stability} becomes larger as the weight coefficient $\tilde{\omega}$ decreases. Therefore, once the stability structure preserving  condition in  (\ref{condition-stability}) holds, the present MRT-LB model (\ref{lb}) must be $L^2$ stable.

    \subsection{The automatic approach to determine the  relaxation parameters and weight coefficients}\label{automatic}
   Based on the above Parts 3.1 and 3.3, it is known that as long as both the fourth-order consistent
   conditions in (\ref{fourth-conition}) and the stability condition in (\ref{condition-stability}) hold, the MRT-LB model (\ref{lb}) for the diagonal-anisotropic diffusion equation (\ref{DE}) can achieve a fourth-order accuracy and is $L^2$ stable. In the following, we will discuss    how to determine the relaxation parameters and weight coefficients from (\ref{fourth-conition}) and (\ref{condition-stability}).

\begin{figure} [H] 
	\begin{center}  
		\subfloat[$\tilde{\omega}=1/18$]
		{
			\includegraphics[width=0.3\textwidth]{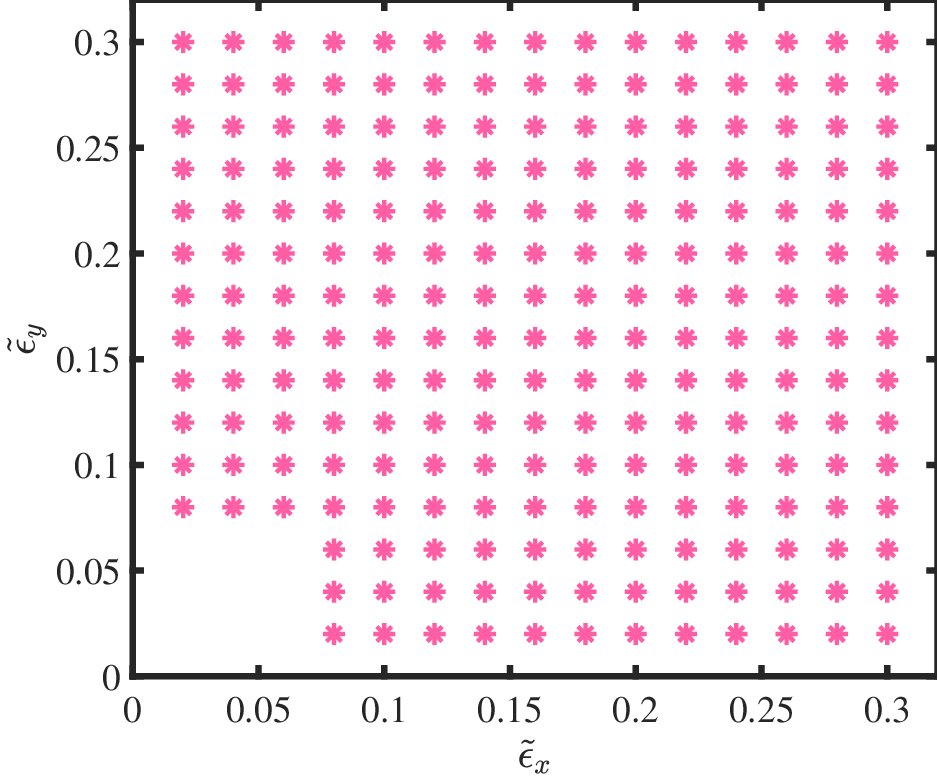}
		} 
		\subfloat[$\tilde{\omega}=1/36$]
		{
			\includegraphics[width=0.3\textwidth]{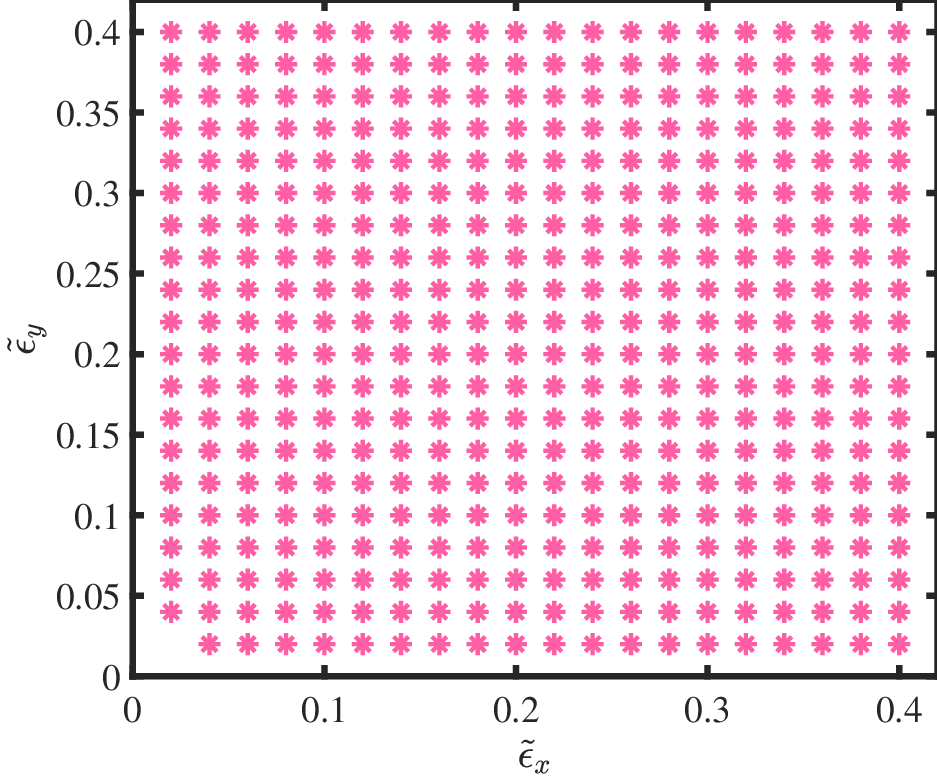}
		}  
		\subfloat[$\tilde{\omega}=1/180$]
		{
			\includegraphics[width=0.302\textwidth]{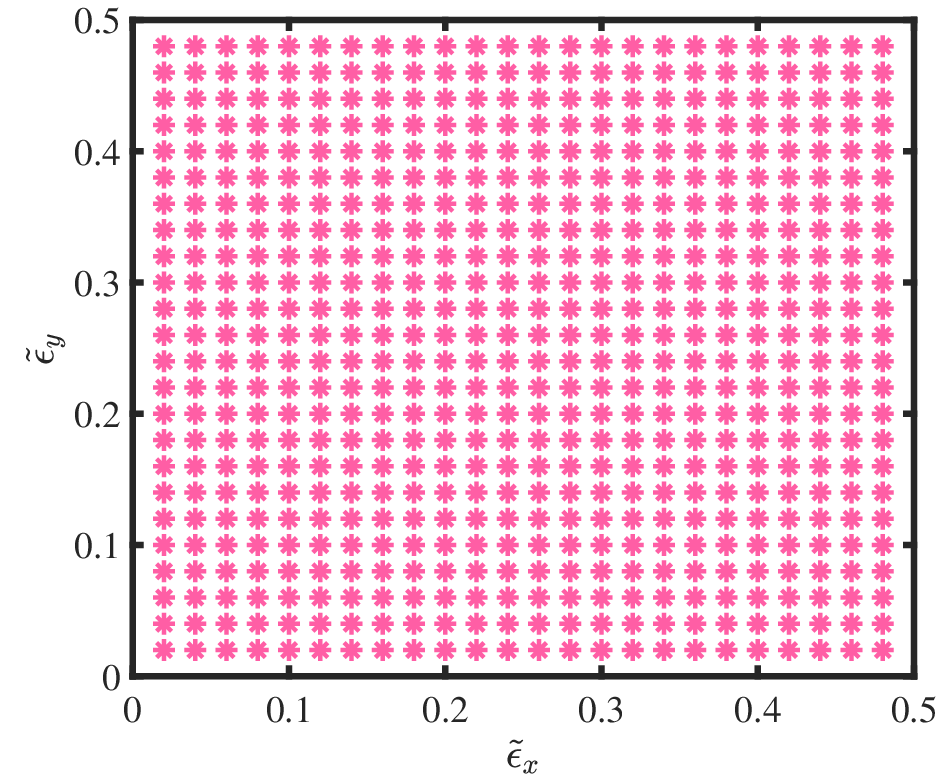}
		}  
		\caption{Regions of the parameter pair $(\tilde{\epsilon}_x,\tilde{\epsilon}_y)$ represent that (\ref{fourth-conition}) and (\ref{condition-stability}) can be solved under  different values of the weight coefficient $\tilde{\omega}$.} \label{fig-region-automatic}\end{center}  
\end{figure}

    According to (\ref{fourth-conition}), one can find that for the general $d$-dimensional case, the choice of $i$ and $j$ is $d(d-1)/2$ in total. In particular, according to the stability condition in  (\ref{condition-stability}), similarly, we can take the relaxation parameter $\tilde{s}_{3|x_i^2x_j}$ in (\ref{s3}) as $\tilde{s}_{x_j}$. Then, in combination with the fourth-order consistent and $L^2$ stability conditions, which can be written as (here we take the choice of $i=j_1$ and $j=j_2$ as an example) 
    \begin{subequations}\label{Gj1j2}
    	\begin{align}   
    		&\tilde{g}_{j_1}(\eta,\Delta t,s_{x_{j_1}},\tilde{s}_{x_{j_1}})=0\quad [\makebox{equivalent to (\ref{fourth-coition-1}) for } i=j_1] ,\label{tgj1}\\ 
    		&\tilde{g}_{j_2}(\eta,\Delta t,s_{x_{j_2}},\tilde{s}_{x_{j_2}})=0\quad [\makebox{equivalent to (\ref{fourth-coition-1}) for } i=j_2],\label{tgj2}\\ &g_{j_1}(\omega_{j_1},\tilde{\omega},\tilde{s}_{x_{j_1}})=0\quad [\makebox{equivalent to (\ref{fourth-coition-2}) for } i=j_1] ,\label{gj1}\\ 
    		&g_{j_2}(\omega_{j_2},\tilde{\omega},\tilde{s}_{x_{j_2}})=0\quad [\makebox{equivalent to (\ref{fourth-coition-2}) for } i=j_2],\label{gj2}\\ 
    		&g_{j_1j_1}(\omega_{j_1},\tilde{\omega},\tilde{s}_{x_{j_1}},s_{2|x_1^2})=0\quad[\makebox{equivalent to (\ref{fourth-coition-3}) and (\ref{condition-stability}) for } i=j_1 ],\label{gj1j1}\\ 
    		&g_{j_2j_2}(\omega_{j_2},\tilde{\omega},\tilde{s}_{x_{j_2}},s_{2|x_1^2})=0\quad[\makebox{equivalent to (\ref{fourth-coition-3}) and (\ref{condition-stability}) for } i=j_2 ],\label{gj2j2}\\ 
    		&g_{j_1j_2}(\omega_{j_1},\omega_{j_2},\tilde{\omega},\tilde{s}_{x_{j_1}},\tilde{s}_{x_{j_2}},s_{2|x_1^2},s_{2|x_{j_1}x_{j_2}})=0\quad[\makebox{equivalent to (\ref{fourth-coition-4}) and (\ref{condition-stability}) for } i=j_1 \:\: j=j_2].\label{gj1j2}
    	\end{align} 
    \end{subequations}
 To solve (\ref{Gj1j2}), in this work, we set the proper weight coefficient  $\tilde{\omega}\in(0,1)$ and relaxation parameter $s_{2|x_1^2}\in(0,2)$ in advance, then $s_{x_{j_1}}$, $s_{x_{j_2}}$, $\omega_{j_1}$, $\omega_{j_2}$, $\tilde{s}_{x_{j_1}}$, $\tilde{s}_{x_{j_2}}$, and $s_{2|x_{j_1}x_{j_2}}$ can be determined by (\ref{tgj1}), (\ref{tgj2}), (\ref{gj1}), (\ref{gj2}), (\ref{gj1j1}), (\ref{gj2j2}), and (\ref{gj1j2}), respectively. In particular, this solving process can be automated by  a simple Matlab symbolic computation code (see \ref{code}), and after running this code $d(d-1)/2$ times (e.g., x-y, x-z, y-z three groups for the case of $d=3$), all the unknown relaxation parameters and weight coefficients  of the MRT-LB model (\ref{lb}) can be determined. To show that the above automatic approach is valid, here we take the two-dimensional case as an example. Specifically, we plot the region of the parameter pair $(\tilde{\epsilon}_{x},\tilde{\epsilon}_{y})$, which represents that  the relaxation parameters and weight coefficients can be obtained from (\ref{fourth-conition}) and (\ref{condition-stability}), in Fig. \ref{fig-region-automatic}, where $s_{2|x_i^2}=1$, $\eta=1$, $\Delta t=1/40$, and different values of the weight coefficient $\tilde{\omega}$ are considered. And from Fig. \ref{fig-region-automatic}, it can be observed that the region becomes larger as the weight coefficient $\tilde{\omega}$ decreases.
Regarding the above discussion on the determination of the relaxation parameters and weight coefficients from (\ref{fourth-conition}) and (\ref{condition-stability}), we now give a remark.\\
\textbf{Remark 3.} As can be seen from (\ref{fourth-coition-1}) and (\ref{fourth-coition-2}), the diffusion coefficient $\kappa_{x_i}$ is associated with the time step $\Delta t$, except when $s_{x_i}=1$.    Thus, for the case of $s_{x_i}=1$, based on the stability condition (\ref{condition-stability}), we further set $s_{3|x_ix_j^2}=1$ and $s_{2|x_i^2}=s_{2|x_1^2}$, (\ref{fourth-conition}) can be further simplified as
\begin{subequations}\label{fourth-coition-sxi-eq0}
	\begin{align}     
				&\omega_i+2(d-1)\tilde{\omega}=\kappa_{x_i}\xi,\quad i\in\llbracket 1,d\rrbracket,\\ 
				&\frac{\omega_i+2(d-1)\tilde{\omega}}{4} =\frac{1}{24}-\frac{1}{4}\Big(1-\frac{1}{s_{2|x_1^2}}\Big)+\Big(\frac{1}{2}-\frac{1}{2s_{2|x_1^2}}\Big)\big(\omega_i+2(d-1)\tilde{\omega}\big),  \quad i\in\llbracket 1,d\rrbracket,\\
				&\Big(\frac{1}{2}+\frac{1}{2s_{2|x_1^2}}\Big)\big(\omega_i+2(d-1)\tilde{\omega}\big)\big(\omega_j+2(d-1)\tilde{\omega}\big) =  \Big(
				\frac{2}{s_{2|x_1^2}}
				+\frac{4}{s_{2|x_ix_j}} -5\Big)\tilde{\omega}  ,\quad i,j\in\llbracket 1,d\rrbracket,\:\: i< j,
			\end{align} 
		\end{subequations} 
		and it is obvious that the above equation has no solution unless  (\ref{DE}) to be solved is reduced to the isotropic type, i.e., the weight coefficient associated with the diffusion coefficient $\kappa_{x_i}$ satisfies $\omega_1=\omega_2=\cdots=\omega_d$. Then,  the relaxation parameters and weight coefficients  can be explicitly obtained from (\ref{fourth-coition-sxi-eq0}) and they are expressed as
		\begin{subequations}\label{si-eq1}
			\begin{align}  	&\omega_0=1-2d\omega_1-2d(d-1)\tilde{\omega},\\
				&\omega_i = \epsilon+2\tilde{\omega}-2d\tilde{\omega},\quad i\in\llbracket 1,d\rrbracket,\\
				&s_{2|x_i^2} = \frac{6(2\epsilon-1)}{6\epsilon-5},\quad i\in\llbracket 1,d\rrbracket,\\
				&s_{2|x_ix_j}= -\frac{6\tilde{\omega}(2\epsilon- 1)(2\epsilon - 1)}{22\epsilon\tilde{\omega} - \epsilon^2 - 5\tilde{\omega}+ 2\epsilon^3- 24\epsilon ^2\tilde{\omega}}\quad i,j\in\llbracket 1,d\rrbracket,\:\: j\neq i,
			\end{align}
		\end{subequations} 
		where $\epsilon\overset{\Delta}{=}\epsilon_{x_i}=\tilde{\epsilon}_{x_i}$ for all $i\in\llbracket 1,d\rrbracket$. This indicates that for the isotropic diffusion equation, the relaxation parameters and weight coefficients, which are used to ensure that the MRT-LB model (\ref{lb}) is fourth-order accurate and $L^2$ stable, are not unique. They can be either automatically  obtained from (\ref{Gj1j2}) without any other restrictions or given by (\ref{si-eq1}) with $s_{x_i}=1$. In addition, we would also like to point out that for the isotropic diffusion equation,  the fourth-order consistent condition (\ref{fourth-coition-ddq2d+1}) of the MRT-LB model with the D$d$Q($2d+1$) lattice structure in Remark 2  can also be explicitly solved, and  the solution is   given by
		\begin{subequations}\label{ddq2d+1-1}
			\begin{align}
				&\omega_i=1-4\sqrt{3}\epsilon,\quad i\in\llbracket 1,d\rrbracket,\\
				&s_{x_i}=\frac{2\Delta t\eta}{\Delta t\eta-\sqrt{(\tilde{s}_{x_i}-4\Delta t\eta+\Delta t^2\tilde{s}_{x_i}\eta^2+2\Delta t\tilde{s}_{x_i}\eta)/\tilde{s}_{x_i}}+1},\quad  i\in\llbracket 1,d\rrbracket,\\
				&s_{2|x_i^2}=4\sqrt{3}-6,  \quad i\in\llbracket 1,d\rrbracket,
			\end{align}
		\end{subequations}
 where $\tilde{s}_{x_i}=6/(3+\sqrt{3}) $ for all $i\in\llbracket 1,d\rrbracket$.

	\section{Numerical results and discussion}\label{Numer}
In this section, we will conduct two benchmarks to  test the numerical accuracy of the present MRT-LB model (\ref{lb}) for the diagonal-anisotropic diffusion equation (\ref{DE}), where the relaxation parameters and weight coefficients  are given by the proposed automatic approach (see Parts \ref{automatic} and \ref{code}), and the numerical accuracy is determined based on the following relative error in $l^2$-norm,
\begin{align*}
	l^2(\phi)=\sqrt{\frac{\sum_i\phi(\mathbf{x}_i,t_n)-\phi^{\star}(\mathbf{x}_j,t_n)}{\sum_i\big(\phi^{\star}(\mathbf{x}_i,t_n)\big)^2)}},
\end{align*}
where $\phi(\mathbf{x}_j,t_n)$ and $\phi^{\star}(\mathbf{x}_j,t_n)$ represent the numerical and analytical solutions at lattice node $\mathbf{x}_j$ ad time $t_n$. 
\subsection{Gauss Hill in two-, three-, and four-dimensional cases}
The first example considers the Gauss Hill problem, and the initial condition is expressed as
\begin{align*}
	\phi(\mathbf{x},0)=\frac{\phi_0}{2\pi|\Gamma_0|}\exp\Big(-\frac{\mathbf{x}^2}{2\Gamma_0^2}\Big).
\end{align*}
With the periodic boundary condition, one can obtain the analytical solution of this problem,
\begin{align*}
	\phi^{\star}(\mathbf{x},t)=\frac{\phi_0}{2\pi|\det(\Gamma)|^{1/2}}\exp\Bigg(-\frac{\sum_{i,j=1}^d(\Gamma^{-1})_{ij}(\mathbf{x}\mathbf{x}^T)_{ij}}{2}\Bigg),
\end{align*}
where $\Gamma=\Gamma_0^2\mathbf{I}+\bm{\kappa}t\in\mathbb{R}^{d\times d} $ [$\bm{\kappa}=\makebox{\textbf{diag}}(\kappa_{x_1},\kappa_{x_2},\cdots,\kappa_{x_d})$ with $d=2,3$, and 4] and  $\det(\Gamma)$ represents the determinant values of $\Gamma$. In the following simulations, the computational domain is set to be $[-1,1]^d$ and the total concentration $\phi_0$ is taken as $\phi_0=2\pi\Gamma_0$, where $\Gamma_0$ should be small enough to ensure that the periodic boundary condition makes sense,  and  we here consider it as 0.05.

\begin{figure} [H] 
	\begin{center}  
		\subfloat[$(\epsilon_{x_1},{\epsilon}_{x_2})=(0.40,0.10)$]
		{
			\includegraphics[width=0.3\textwidth]{ 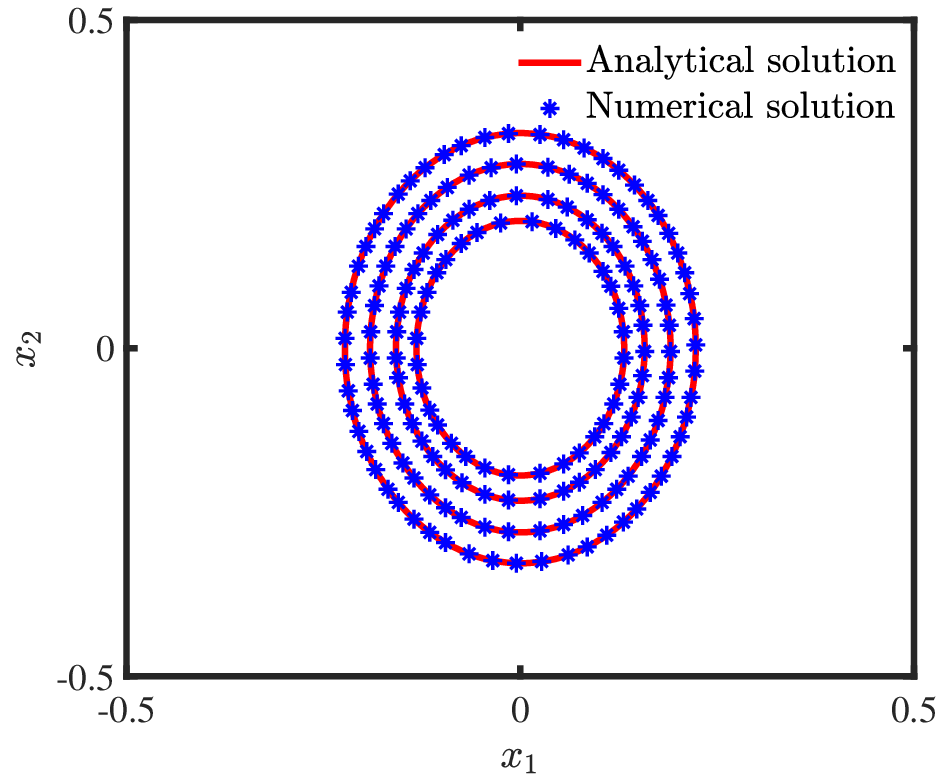}
		} 
		\subfloat[$(\epsilon_{x_1},{\epsilon}_{x_2})=(0.20,0.20)$]
		{
			\includegraphics[width=0.3\textwidth]{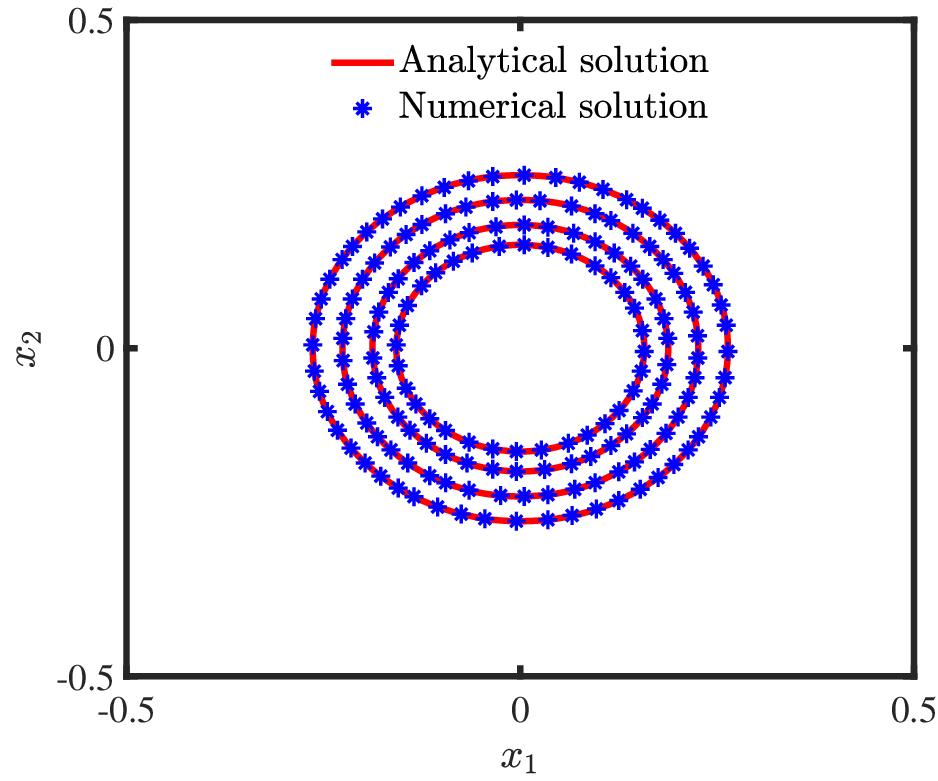}
		}  
		\subfloat[$(\epsilon_{x_1},{\epsilon}_{x_2})=(0.10,0.30)$]
		{
			\includegraphics[width=0.3\textwidth]{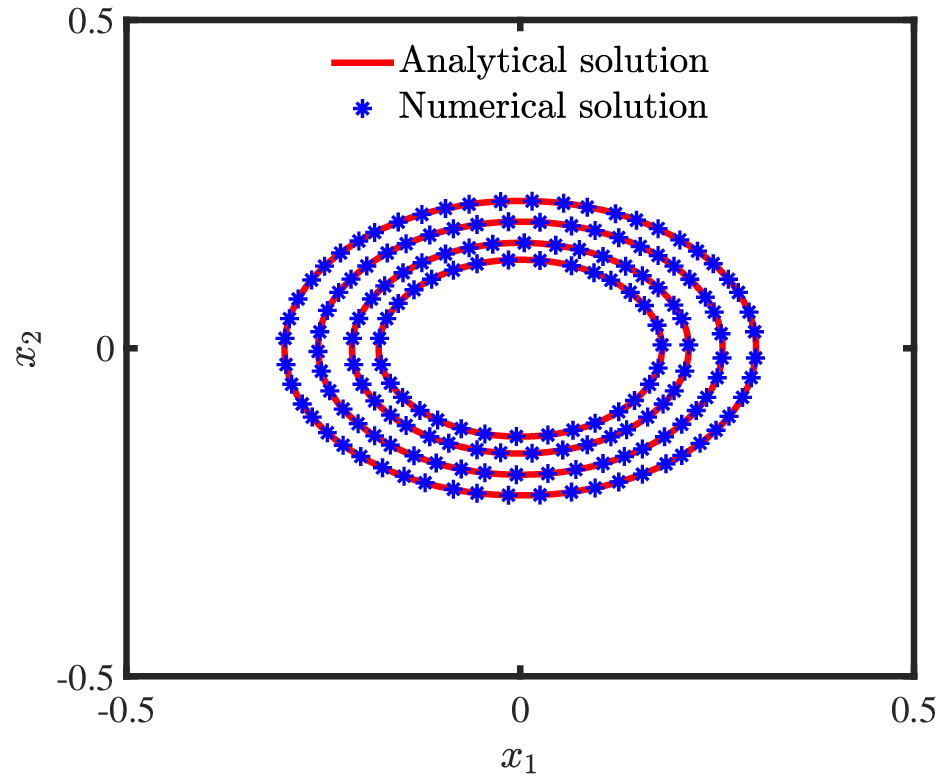}
		}  
		\caption{ Contour lines of the numerical and analytical solutions under different values of the parameter pair $(\epsilon_{x_1},{\epsilon}_{x_2})$.} \label{fig-ex1}\end{center}  
\end{figure}

In order to show the performance of the present MRT-LB model, we first conduct some tests with different values of the diagonal diffusion matrix $\bm{\kappa}$ in two-dimensional case,  where the lattice spacing $\Delta x=1/100$, the  time step $\Delta t=1/40$, and the weight coefficients and   relaxation parameters  are presented in Table \ref{table-para-ex1-aniso}. The contour lines of the numerical and analytical solutions at the final time $t=2$ are plotted in Fig. \ref{fig-ex1}, from which one can find that the numerical results of the MRT-LB model are in good agreement with the analytical solutions.

\begin{figure} [H] 
 	\begin{center}  
 		\subfloat[$d=2(\Delta x^2/\Delta t=1/250)$]
 		{
 			\includegraphics[width=0.3\textwidth]{ 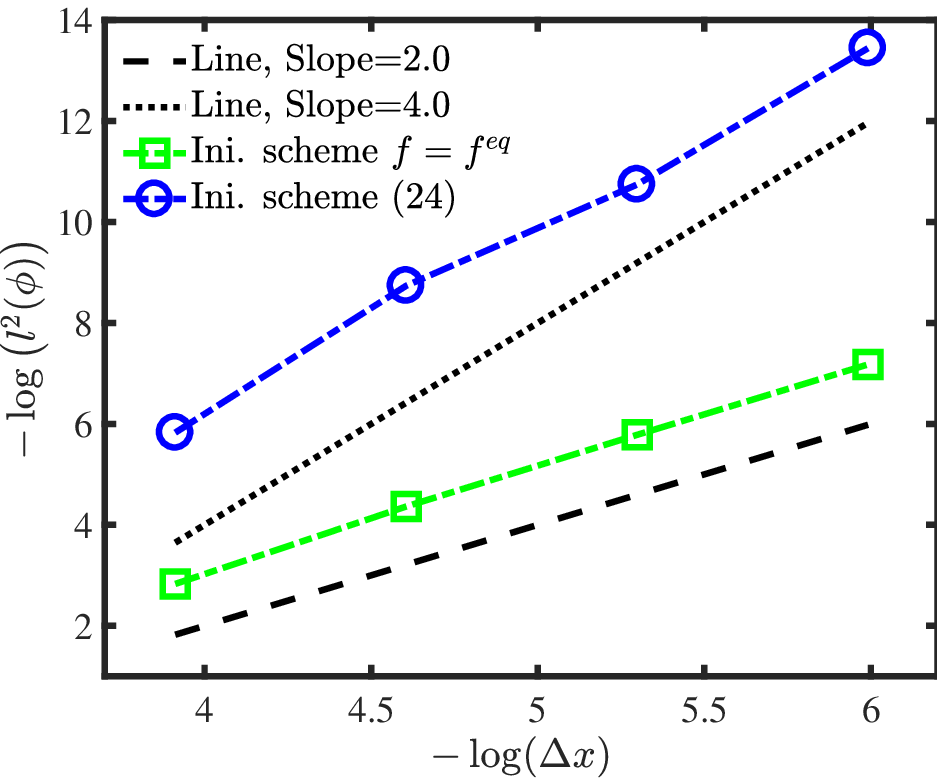}
 		} 
 		\subfloat[$d=3(\Delta x^2/\Delta t=1/40)$]
 		{
 			\includegraphics[width=0.3\textwidth]{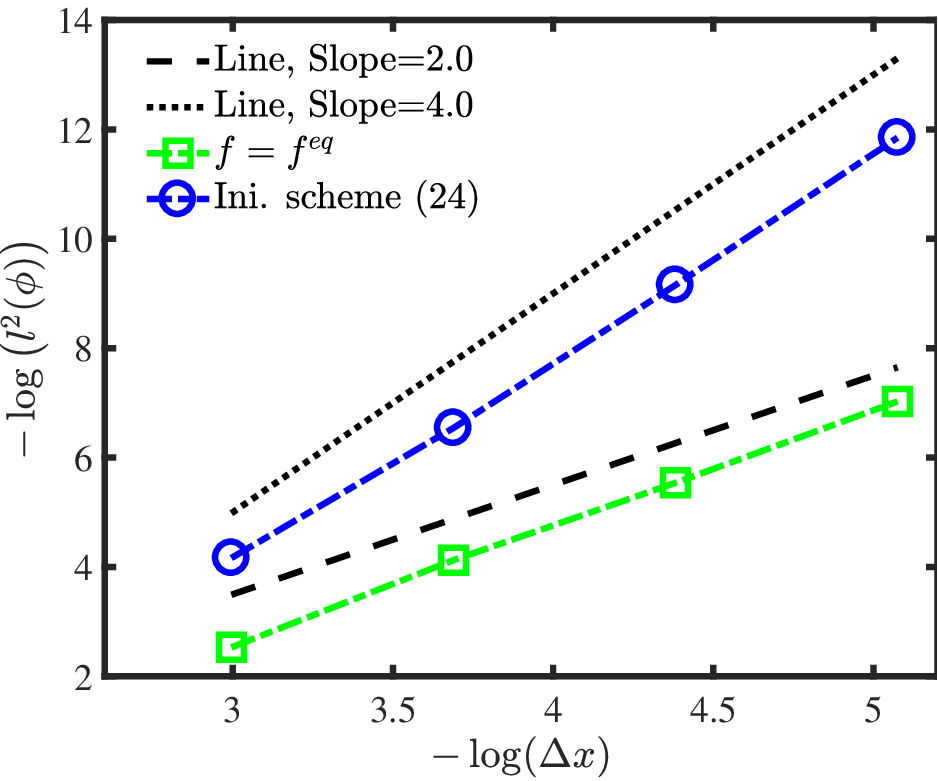}
 		}  
 		\subfloat[$d=4(\Delta x^2/\Delta t=1/50)$]
 		{
 			\includegraphics[width=0.3\textwidth]{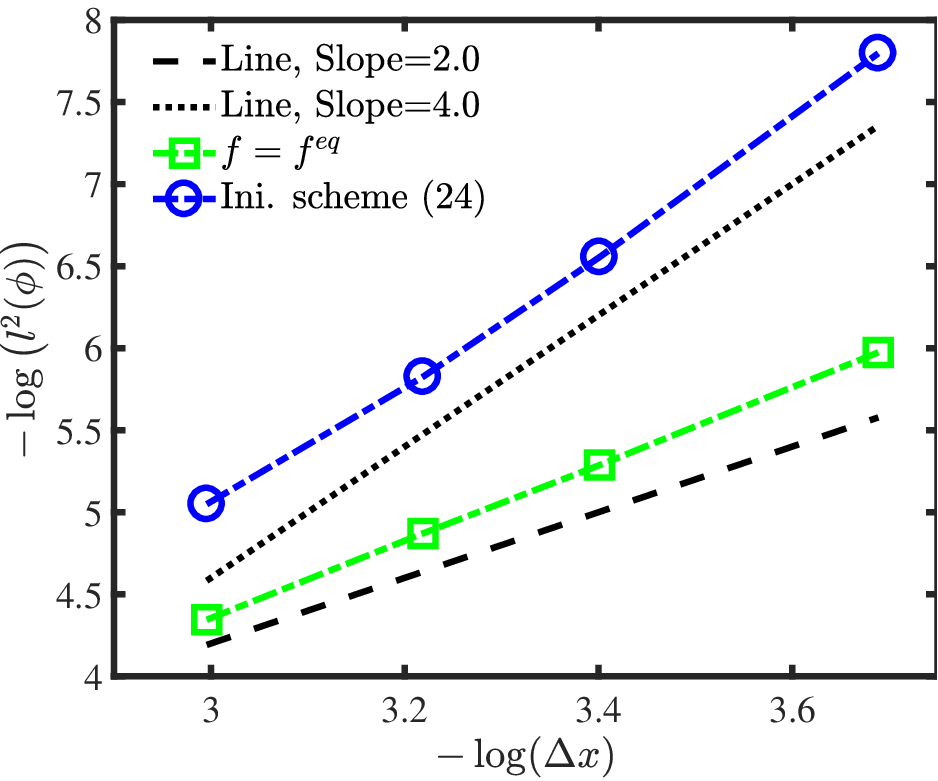}
 		}  
 		\caption{Convergence rates of the MRT-LB model with different initialization schemes at the final time $t=2$.} \label{fig-ex1-slope}\end{center}  
 \end{figure}

Then,  to investigate the impact of the initialization scheme on the accuracy of the MRT-LB model, for the Gauss Hill in all the two-, three-, and four-dimensional cases, we test the convergence rates of the present MRT-LB model with different initialization schemes of $f_k=f_k^{eq}$ and (\ref{ini}). The results are presented in Fig. \ref{fig-ex1-slope}, where $(\epsilon_{x_1},\epsilon_{x_2})=(0.30,0.10)$, $(\epsilon_{x_1},\epsilon_{x_2},\epsilon_{x_3})=(0.10,0.40,0.15)$, $(\epsilon_{x_1},\epsilon_{x_2},\epsilon_{x_3},\epsilon_{x_4})=(0.15,0.20,0.10,0.05)$, and the values of  relaxation parameters and weight coefficients for $d=3$ and 4 are provided in Table \ref{table-para-ex1-aniso}. From this figure, we can clearly see that when the distribution function is initialized at the equilibrium state $f_k^{eq}$, the MRT-LB model only achieves a second-order accuracy. However,  when the initialization scheme in (\ref{ini}) is adopted, the MRT-LB model can be fourth-order accurate, which is in full agreement with our theoretical analysis in Part \ref{ini-sec}.

\begin{table}[h]   
	\small
	\begin{center}   
		\caption{ Relaxation parameters and weight coefficients of the MRT-LB  model for $d=2$, 3, and 4.}  	 
		\label{table-para-ex1-aniso}
		\begin{tabular}{cccccccccc}\hline\hline  
			$({\epsilon}_{x_1},{\epsilon}_{x_2},\cdots,{\epsilon}_{x_d})$&$(0.40,0.10)$&$(0.20,0.20)$&$(0.10,0.30)$&$(0.10,0.40,0.15)$&$(0.15,0.20,0.10,0.05)$ \\
			\hline
			$\omega_0$&0.392414074930637&0.466081465800228 & 0.279164263737521 &0.044310556197977 &0.003107936020711\\
			$\omega_1$&0.003792962534682 & 0.105701855772165&11/45    &1/9 &0.148170462855893 \\
			$\omega_2$&11/45&0.105701855772165&0.060417868131240 &0.037126295868015&0.144590744661054 \\
			$\omega_3$&-&-&-&0.296273981588552&0.116666666666667 \\
			$\omega_4$&-&-&-&-&0.055684824472697\\
			$\tilde{\omega}$&1/36&1/36&1/36&1/180&1/360\\
			$s_{x_1}$&0.258403002308493&0.892756279989137&3/2 & 8/7& 1.047126365130629 \\
			$s_{x_2}$&3/2&0.892756279989137& 0.557600159447285 &0.258403002308493&   0.892756279989137\\
			$s_{x_3}$&-&-&-&1.359653295886320 &1.142857142857143  \\
			$s_{x_4}$&-&-&-&-& 1.182682621447616  \\
			$s_{2|x_1^2}$&1&1 &1&1&1 \\
			$s_{2|x_1x_2}$&1.466835061000191&1.011385583524664&1.192683097984767  & 0.945790034643835 &0.299130236472667\\
			$s_{2|x_1x_3}$&-&-&-&1.151202850452001 &0.485974551112802\\
			$s_{2|x_1x_4}$&-&-&-&- &0.696896856214742\\
			$s_{2|x_2x_3}$&-&-&-&0.770241927190338 &0.408239754101923\\
			$s_{2|x_2x_4}$&-&-&-&- &0.625878745350766\\
			$s_{2|x_3x_4}$&-&-&-&- &0.812554973056151\\
			\hline\hline
		\end{tabular}
	\end{center}  
\end{table} 

\subsection{Diagonal-anisotropic diffusion equation with a linear source}
In the second example,  we focus on the two-dimensional diagonal-anisotropic diffusion equation (\ref{DE}) with a linear source term, in which the  initial condition is given by
\begin{align*}
	\phi(x,y,0)=\sin(\pi x)\sin(\pi y)+\pi^2,\quad (x,y)\in(-1,1)\times (-1,1),
\end{align*}
and the analytical solution can be obtained as
\begin{align*}
	\phi^{\star}(x,y,t)=\sin(\pi x)\sin(\pi y)\exp\Big[\big(-\pi^2(\kappa_x+\kappa_y+1)\big)\Big].
\end{align*}

\begin{figure} [H] 
	\begin{center}  
		\subfloat[$(\tilde{\epsilon}_x,\tilde{\epsilon}_y)=(0.25,0.10)$]
		{
			\includegraphics[width=0.3\textwidth]{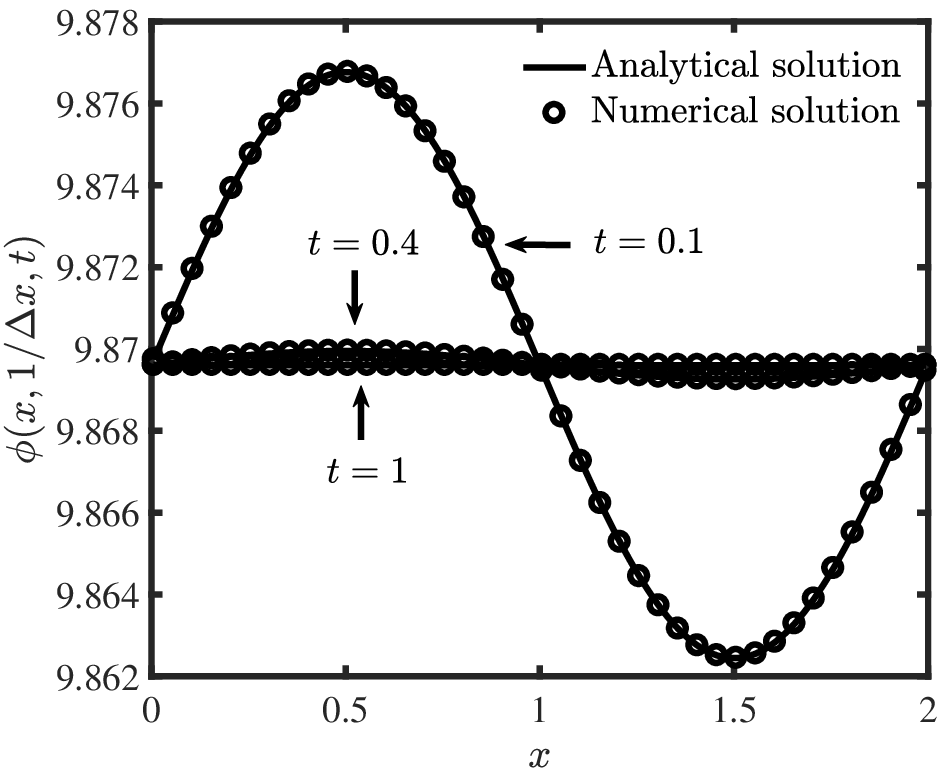}
		} 
		\subfloat[$(\tilde{\epsilon}_x,\tilde{\epsilon}_y)=(0.40,0.40)$]
		{
			\includegraphics[width=0.3\textwidth]{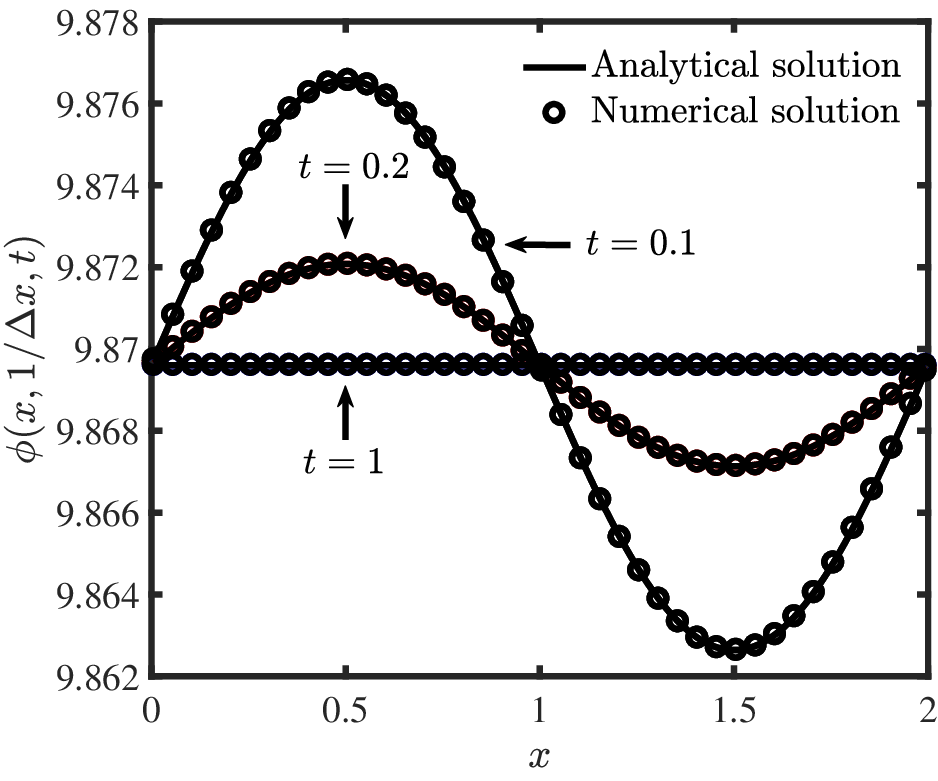}
		}  
			\subfloat[$(\tilde{\epsilon}_x,\tilde{\epsilon}_y)=(0.15,0.40)$]
		{
			\includegraphics[width=0.3\textwidth]{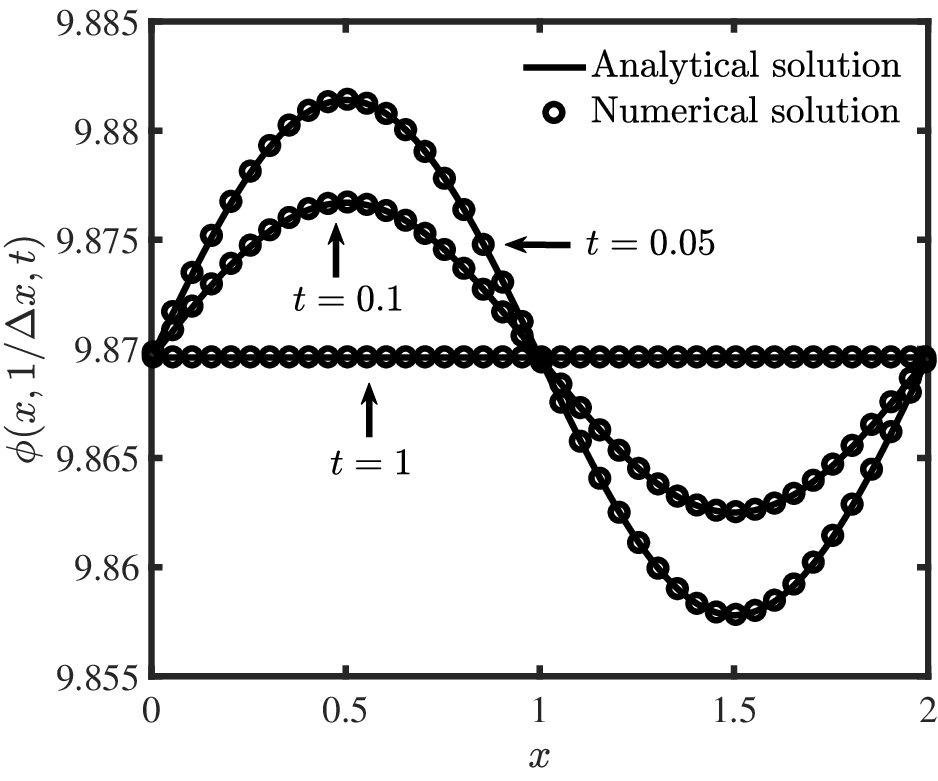}
		}  
		\caption{Profiles of the numerical and analytical solutions under different values of the final time $t$.} \label{fig-ex2}\end{center}  
\end{figure}

\begin{figure} [H] 
 	\begin{center}  
 			\includegraphics[width=0.5\textwidth]{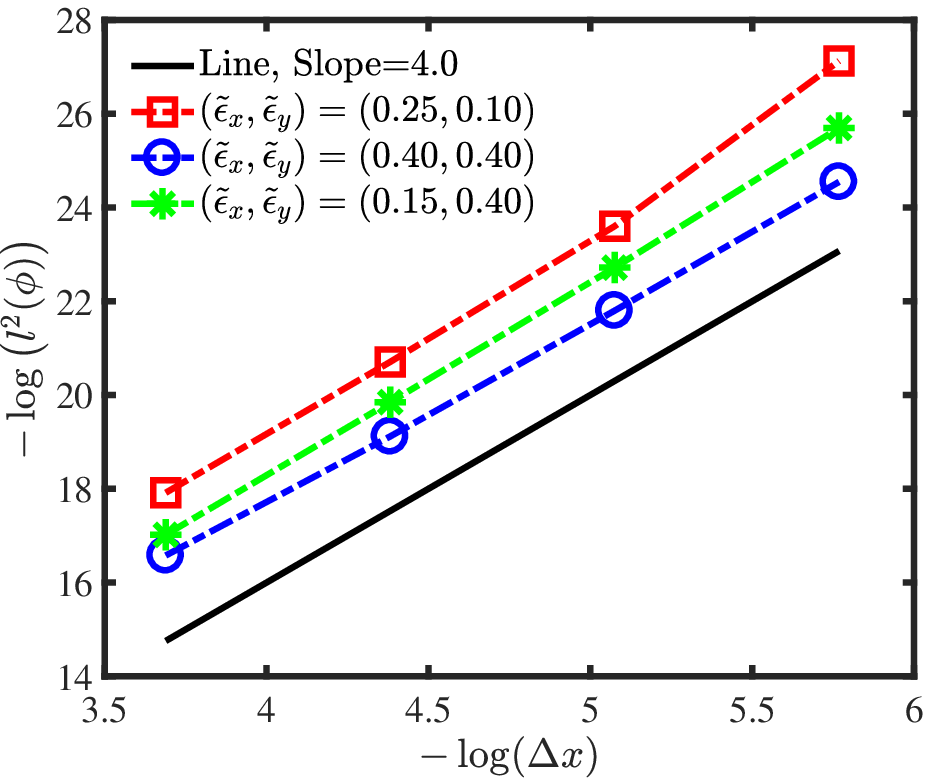}
 		\caption{Convergence rates of the MRT-LB model under different values of the parameter pair ($\tilde{\epsilon}_x,\tilde{\epsilon}_y$).} \label{fig-ex2-slope}\end{center}  
 \end{figure}

First of all, we consider three different values of the parameter pair $(\tilde{\epsilon}_x,\tilde{\epsilon}_y)=(0.25,0.10)$, $(0.40,0.40)$, and $(0.15,0.40)$. For the MRT-LB model, the weight coefficient $\tilde{\omega}=1/36$, the relaxation parameter $s_{2|x_i^2}=1$, and the other weight coefficients and relaxation parameters can be seen in Table \ref{table-para-ex2-aniso}. Fig. \ref{fig-ex2} plots the profiles of the numerical results at different values of the final time $t$, where the lattice spacing $\Delta x=1/80$ and the time step $\Delta t=1/400$. From Fig. \ref{fig-ex2}, one can see that the numerical results are in good agreement with the analytical solutions. In addition, the logarithm of the relative errors at the final time $t=1$ with respect to the lattice spacing $\Delta x$ is also shown in Fig. \ref{fig-ex2-slope}, where the time step $\Delta t$ is given by the relation $\Delta x^2/\Delta t=1/16$. From this figure, it can be seen that for the diagonal-anisotropic diffusion equation with a linear source term, the present MRT-LB model can also be fourth-order accurate, which is consistent with our theoretical analysis. 
\begin{table}[h]   
	\begin{center}   
		\caption{ Relaxation parameters and weight coefficients of the MRT-LB  model under different values of the parameter pair $(\tilde{\epsilon}_x,\tilde{\epsilon}_y)$.}  	 
		\label{table-para-ex2-aniso}
		\begin{tabular}{cccccccccc}\hline\hline  
				$(\tilde{\epsilon}_x,\tilde{\epsilon}_y)$&$(0.25,0.10)$&$(0.40,0.40)$&$(0.15,0.40)$ \\
			\hline
			$\omega_0$&    0.228240384066444
						&0.873717038750163  &    0.662739815885518
						\\
			$\omega_1$&   0.085879807966778
			&0.003792962534682 &0.109281573967004 \\
			$\omega_2$&11/45&0.003792962534682&   0.003792962534682
			\\
			$s_{x}$&   0.729269281934827
			&0.314095759114162 &    1.045990365920910
			 \\
			$s_{y}$&   1.487864247860460
			&0.314095759114162&   0.275195555816491
			\\
		    $s_{2|xy}$&   1.114757496216724
		    &1.770421857439279  &   1.468455215964528
		    \\
			\hline\hline
		\end{tabular}
	\end{center}  
\end{table}

Then, when the diagonal-anisotropic diffusion equation is reduced to the isotropic type, i.e., $\tilde{\epsilon}_x=\tilde{\epsilon}_y$, a simulation is conducted to give a comparison of the fourth-order MRT-LB-D2Q9-(\ref{Gj1j2}) \footnotemark[2],   MRT-LB-D2Q9-(\ref{si-eq1}) \footnotemark[3], and  MRT-LB-D2Q5-(\ref{ddq2d+1-1}) \footnotemark[4] models. Considering two different values of the parameter $\tilde{\epsilon}_x=0.1$ and $0.2$, we measure the relative errors of these three models and calculate the convergence rates (CRs) in Table \ref{table-iso-ex2-1}, where the final time $t=1$ and the weight coefficients and relaxation parameters of these three models are provided in Table \ref{table-para-ex2}. From Table \ref{table-iso-ex2-1} with the parameter  $\tilde{\epsilon}_x=0.1$, one can find that all the three MRT-LB models have a fourth-order accuracy, which is in agreement with our theoretical analysis in Remark 3. However, for the   case of $\tilde{\epsilon}_x=0.2$ shown in Table \ref{table-iso-ex2-1}, both the MRT-LB-D2Q9-(\ref{Gj1j2}) and MRT-LB-D2Q9-(\ref{si-eq1}) models can achieve a fourth-order accuracy while the MRT-LB-D2Q5-(\ref{ddq2d+1-1}) model cannot.  This is due to the fact that, as can be seen in Table \ref{table-para-ex2}, the weight coefficient $\omega_0<0$ of the MRT-LB-D2Q5-(\ref{ddq2d+1-1}) model, which means that the MRT-LB-D2Q5-(\ref{ddq2d+1-1}) model is no longer stable for the case of  $\tilde{\epsilon}_x=0.2$. Therefore, we can conclude that for the isotropic diffusion equation, the fourth-order  MRT-LB  model with the D$d$Q($2d^2+1$) lattice structure  is more general  than that with the D$d$Q($2d+1$) lattice structure .

\begin{table}[h] 
	\small  
	\begin{center}   
		\caption{ Relaxation parameters and weight coefficients of the MRT-LB-D2Q9-(\ref{Gj1j2}), MRT-LB-D2Q9-(\ref{si-eq1}), and MRT-LB-D2Q5-(\ref{ddq2d+1-1}) models.}  	 
		\label{table-para-ex2}
		\begin{tabular}{cccccccccc}\hline\hline  
			Model&\multicolumn{2}{c}{MRT-LB-D2Q9-(\ref{Gj1j2})}&&\multicolumn{2}{c}{MRT-LB-D2Q9-(\ref{si-eq1})}&&\multicolumn{2}{c}{MRT-LB-D2Q5-(\ref{ddq2d+1-1})} \\
			\cline{2-3}\cline{5-6}\cline{8-9}
			$\tilde{\epsilon}_x$&0.1&0.2&&0.1&0.2&&0.1&0.2\\
			\hline
			$\omega_0$&4/9&4/9&&32/45&14/45&&$1-0.4\sqrt{3}$&$(1-0.8\sqrt{3})<0$\\
			$\omega_1=\omega_2$&1/9&1/9&&2/45&13/90 &&$0.1\sqrt{3}$&$0.2\sqrt{3}$\\
			$\tilde{\omega}$&1/36&1/36&&1/36&1/36&&-&-\\
			$s_{x_1}=s_{x_2}$&1.227162705974765&0.917191588597811&&1&1&&1.243448367507882&1.261464647585203\\
			$s_{2|x_1^2}=s_{2|x_2^2}$&1&1&&12/11&18/19 &&$4\sqrt{3}-6$&$4\sqrt{3}-6$\\
			$s_{2|x_1x_2}$&0.948397105014222 &0.993142767001046 &&15/13&90/101&&-&-\\
			\hline\hline
		\end{tabular}
	\end{center}  
\end{table}

\begin{table}[h]   
	\begin{center}   
		\caption{ Relative errors and convergence rates of the MRT-LB-D2Q9-(\ref{Gj1j2}), MRT-LB-D2Q9-(\ref{si-eq1}), and MRT-LB-D2Q5-(\ref{ddq2d+1-1}) models under different values of the parameter $\tilde{\epsilon}_x$.}  	 
		\label{table-iso-ex2-1}
		\begin{tabular}{cccccccccc}\hline\hline  
			\multirow{2}{*}{$\tilde{\epsilon}_x$}&\multirow{2}{*}{$\Delta x$}&\multicolumn{2}{c}{MRT-LB-D2Q9-(\ref{Gj1j2})}&&\multicolumn{2}{c}{MRT-LB-D2Q9-(\ref{si-eq1})}&&\multicolumn{2}{c}{MRT-LB-D2Q5-(\ref{ddq2d+1-1})} \\
			\cline{3-4}\cline{6-7}\cline{9-10}
			&&$l^2(\phi)$&CR&&$l^2(\phi)$&CR&&$l^2(\phi)$&CR\\
			\hline
			\multirow{5}{*}{0.1}&1/40&1.8240$\times 10^{-8}$&-&&1.8101$\times 10^{-8}$&-&&1.8539$\times 10^{-8}$&-\\
			&1/80&1.1344$\times 10^{-9}$&3.9961&&1.1614$\times 10^{-9}$&3.9966 &&1.1318$\times 10^{-9}$&3.9958\\
			&1/160&7.0916$\times 10^{-11}$&3.9996 &&7.2597$\times 10^{-11}$&3.9998&&7.0760$\times 10^{-11}$&3.9996\\
			&1/240&1.4008$\times 10^{-11}$&3.9999 &&1.4340$\times 10^{-11}$&4.0000&&1.3978$\times 10^{-11}$&3.9999\\
			&1/320&4.4324$\times 10^{-12}$&4.0000 &&4.5375$\times 10^{-12}$&3.9999&&4.4227$\times 10^{-12}$&4.0000\\
			\hline
			\multirow{5}{*}{0.2}&1/40&1.6813$\times 10^{-8}$&-&&1.6331$\times 10^{-8}$&-&&0.9895&-\\
			&1/80&1.0549$\times 10^{-9}$&3.9944 &&1.0232$\times 10^{-9}$& 3.9965&&1.0000&-0.0152\\
			&1/160&6.5971$\times 10^{-11}$&3.9992 &&6.3962$\times 10^{-11}$&3.9998&&-&-\\
			&1/240&1.3032$\times 10^{-11}$&3.9998 &&1.2634$\times 10^{-11}$&4.0000&&-&-\\
			&1/320&4.1237$\times 10^{-12}$&3.9999 &&3.9976$\times 10^{-12}$&4.0000&&-&-\\
			\hline\hline
		\end{tabular}
	\end{center}  
\end{table}

	\footnotetext[2]{The MRT-LB model with the D2Q9 lattice structure and the weight coefficients and relaxation parameters automatically determined  by (\ref{Gj1j2}).}
	\footnotetext[3]{The MRT-LB model with the D2Q9 lattice structure and  the weight coefficients and relaxation parameters given  by (\ref{si-eq1}).}
\footnotetext[4]{The MRT-LB model with the D2Q5 lattice structure and the weight coefficients and relaxation parameters given by (\ref{ddq2d+1-1}).}

	\section{Conclusion}\label{conclusion}
  This work proposes an automatic approach to determine the  fourth-order and $L^2$ stable MRT-LB model for the diagonal-anisotropic diffusion equation,  in which the MRT-LB model adopts the transformation matrix constructed in a natural way and the D$d$Q($2d^2+1$) lattice structure. 
  \begin{itemize}
  	\item Through the direct Taylor expansion, the conditions, which are associated with the relaxation parameters and weight coefficients, are derived to ensure that the macroscopic modified equation of the MRT-LB model can be fourth-order consistent with the diagonal-anisotropic diffusion equation. In particular, for the isotropic diffusion equation, another MRT-LB model with the D$d$Q($2d+1$) lattice structure [fewer discrete velocities than the D$d$Q($2d^2+1$) lattice structure] is developed, and the fourth-order consistent conditions are also similarly given. Furthermore, to obtain an overall fourth-order MRT-LB model, the fourth-order   initialization scheme of the distribution function is also constructed. 
  	\item The condition which guarantees that the present MRT-LB model can satisfy the stability structure is given. In particular, once the stability structure is satisfied, it can be numerically  validated that the MRT-LB model is $L^2$ stable. Then, based on this stability structure preserving condition and the above fourth-order consistent conditions, the relaxation parameters and weight coefficients of the MRT-LB model can be automatically determined  by a simple computer code.
  	\item Some numerical experiments are carried out to test the convergence rate of the proposed MRT-LB model, and the numerical results are consistent with our
  	theoretical analysis.  In addition,  a comparison simulation is conducted with respect to  the isotropic diffusion equation, and the numerical results show that compared with the fourth-order MRT-LB model with the D$d$Q($2d+1$) lattice structure,  the present fourth-order MRT-LB model with the D$d$Q($2d^2+1$) lattice structure is more general.
  \end{itemize}  
  In summary, for the general diagonal-anisotropic diffusion equation, after obtaining the relaxation parameters and weight coefficients by a simple computer code, the MRT-LB model is fourth-order accurate, $L^2$ stable, and also valid for the high-dimensional case.
 
 \appendix
 \section{The second- to fifth-order expansions of the distribution function}\label{sec-app-1}
 To obtain the second- to fifth-order expansions of the distribution function $f_k$,  we here divide the derivation into four steps.
 \subsection{The second-order expansion}
 Substituting  $f_k^{ne}=O(\Delta x)$ into (\ref{f-o6}) yields the following second-order expansion of the distribution function $f_k$, 
 	\begin{align}\label{2rd}
 		 O(\Delta x^2) +\Delta x\mathbf{e}_k\cdot\nabla f_k^{eq}
 		= 	-\sum_{i=1}^q\mathbf{\Lambda}_{ki}f_i^{ne}.
 	\end{align}  
 \subsection{The third-order expansion}
Based on the above (\ref{2rd}), one can obtain the third-order expansion of the distribution function $f_k$ from (\ref{f-o6}),  
 	\begin{align} \label{3th}
 		&\xi\Delta x^2\partial_tf_k^{eq}+\Delta x\mathbf{e}_k\cdot\nabla f_k^{eq}-\Delta x^2\nabla^2\overset{2}{\cdot}\mathbf{e}_k\sum_{i=1}^q\bm{\Lambda}^{-1}_{ki}\mathbf{e}_if_i^{eq}+\frac{\Delta x^2}{2}\nabla^2\overset{2}{\cdot}\mathbf{e}_k^2 f_k^{eq}\notag\\
 		&-\xi\Delta x^2\sum_{i=1}^q\Big(\mathbf{I}-\frac{\mathbf{\Lambda}}{2}\Big)_{ki}R_i+O(\Delta x^3) =  -\sum_{i=1}^q\mathbf{\Lambda}_{ki}f_i^{ne}.
 	\end{align}  
 \subsection{The fourth-order expansion}
 From (\ref{f-o6}), the fourth-order expansion of the distribution function $f_k$ is given by 
  \begin{align} \label{4th}
  	&\xi\Delta x^2\partial_tf_k^{eq}-\xi\Delta x^3\sum_{i=1}^q\mathbf{\Lambda}_{ki}^{-1}\mathbf{e}_i\cdot\nabla\partial_tf_i^{eq} +O(\Delta x^4)=-\Delta x\mathbf{e}_k\cdot\nabla f_k^{eq}+\frac{\Delta x^2}{2}\mathbf{e}_k^2\overset{2}{\cdot}\nabla^2 f_k^{eq}  -\Delta x^2\nabla^2\overset{2}{\cdot}\mathbf{e}_k \sum_{i=1}^q\Big(\mathbf{I}-\bm{\Lambda}^{-1}\Big)_{ki}\mathbf{e}_if_i^{eq}\notag\\
  	&-\sum_{i=1}^q\mathbf{\Lambda}_{ki}f_i^{ne} -\xi\Delta x^3\nabla\cdot\mathbf{e}_k\sum_{i=1}^q\Big(\mathbf{I}-\bm{\Lambda}^{-1}\Big)_{ki}\partial_tf_i^{eq}\notag\\
  	&-\Delta x\nabla^3\overset{3}{\cdot}\mathbf{e}_k\sum_{i=1}^q\bigg(\bm{\Lambda}^{-1}_{ki} \mathbf{e}_i\sum_{l=1}^q\Big(\bm{\Lambda}^{-1}_{il}\mathbf{e}_lf_l^{eq}\Big)\bigg)
  	+\frac{\Delta x^3}{2}\nabla^3\overset{3}{\cdot}\mathbf{e}_k\sum_{i=1}^q\bm{\Lambda}^{-1}_{ki}\mathbf{e}^2_if_i^{eq}
  	+\frac{\Delta x^3}{2}\nabla^3\overset{3}{\cdot}\mathbf{e}_k^2\sum_{i=1}^q\bm{\Lambda}^{-1}_{ki}\mathbf{e}_if_i^{eq}
  	-\frac{\Delta x^3}{6}\nabla^3\overset{3}{\cdot}\mathbf{e}_k^3f_k^{eq} \notag\\
  	&+\xi\Delta x^2\sum_{i=1}^q\Big(\mathbf{I}-\frac{\mathbf{\Lambda}}{2}\Big)_{ki}R_i
  	-\xi\Delta x^3\mathbf{e}_k\cdot\nabla\sum_{i=1}^q\Big(\mathbf{\Lambda}^{-1}-\frac{\mathbf{I}}{2}\Big)_{ki}R_i ,
  \end{align}
 where  (\ref{2rd})  and (\ref{3th}) have been used. After taking the second-order moment  of the above equation, one can obtain the unknown term $\sum_{k=1}^q\mathbf{e}_{ki}\mathbf{e}_{kj}f_k^{ne}$ [see (\ref{fneq-4rd})].
 \subsection{The fifth-order expansion} 
 According to (\ref{2rd}), (\ref{3th}), and (\ref{4th}), we can further obtain the fifth-order expansion of the distribution function $f_k$ from (\ref{f-o6}),
  \begin{align}\label{5th}
 	& \xi\Delta x^2\partial_tf_k^{eq}-\xi\Delta x^2\partial_t\sum_{i=1}^q\bm{\Lambda}^{-1}_{ki}\Big(\xi\Delta x^2\partial_tf_i^{eq}+\Delta x\mathbf{e}_i\cdot\nabla f_i^{eq}-\Delta x^2\nabla^2\overset{2}{\cdot}\mathbf{e}_i\sum_{l=1}^q\bm{\Lambda}^{-1}_{il}\mathbf{e}_lf_l^{eq}+\frac{\Delta x^2}{2}\nabla^2\overset{2}{\cdot}\mathbf{e}_i^2 f_i^{eq}\notag\\
 	&\qquad\qquad\qquad\qquad\qquad\quad-\xi\Delta x^2\sum_{i=1}^q\Big(\mathbf{I}-\frac{\mathbf{\Lambda}}{2}\Big)_{ki}R_i\Big)+\frac{\xi^2\Delta x^4}{2}\partial_t^2f_k^{eq}+O(\Delta x^5)\notag\\
 	&= -\Delta x\mathbf{e}_k\cdot\nabla f_k^{eq}+\frac{\Delta x^2}{2}\mathbf{e}_k^2\overset{2}{\cdot}\nabla^2 f_k^{eq}-\frac{\Delta x^3}{6}\mathbf{e}_k^3\overset{3}{\cdot}\nabla^3 f_k^{eq}
 	+\frac{\Delta x^4}{24}\mathbf{e}_k^4\overset{4}{\cdot}\nabla^4 f_k^{eq} -\sum_{i=1}^q\mathbf{\Lambda}_{ki}f_i^{ne}\notag\\
 	&-\Delta x\mathbf{e}_k\cdot\nabla\sum_{i=1}^q\Big( \mathbf{I}-\mathbf{\Lambda}^{-1}\Big)_{ki}\Bigg(\xi\Delta x^2\partial_tf_i^{eq}-\xi\Delta x^3\sum_{l=1}^q\bm{\Lambda}^{-1}_{il}\mathbf{e}_l\cdot\nabla\partial_tf_l^{eq}+\Delta x\mathbf{e}_i\cdot\nabla f_i^{eq}  -\frac{\Delta x^2}{2}\mathbf{e}_i^2\overset{2}{\cdot}\nabla^2 f_i^{eq} \notag\\
 	&\qquad\qquad\qquad\qquad\qquad\quad+\Delta x^2\nabla^2\overset{2}{\cdot}\mathbf{e}_i\sum_{l=1}^q\Big(\mathbf{I}-\mathbf{\Lambda}^{-1}\Big)_{il}\mathbf{e}_l f_l^{eq}+\xi\Delta x^3\nabla\cdot\mathbf{e}_i\sum_{l=1}^q\Big(\mathbf{I}-\mathbf{\Lambda}^{-1}\Big)_{il}  \partial_tf_l^{eq}
 	\notag\\
 	&\qquad\qquad\qquad\qquad\qquad\quad
 	+\Delta x^3\nabla^3\overset{3}{\cdot}\mathbf{e}_i\sum_{l=1}^q\mathbf{\Lambda}^{-1}_{il} \mathbf{e}_l\sum_{o=1}^q\mathbf{\Lambda}^{-1}_{lo}\mathbf{e}_o f_o^{eq} 
 	-\frac{\Delta x^3}{2}\nabla^3\overset{3}{\cdot}\mathbf{e}_i\sum_{l=1}^q\mathbf{\Lambda}^{-1}_{il}  \mathbf{e}_l^2f_l^{eq}-\frac{\Delta x^3}{2}\nabla^3\overset{3}{\cdot}\mathbf{e}_i^2 \sum_{l=1}^q\mathbf{\Lambda}_{il}^{-1}\mathbf{e}_l f_l^{eq} 
 	\notag\\
 	&\qquad\qquad\qquad\qquad\qquad\quad+\frac{\Delta x^3}{6}\nabla^3\overset{3}{\cdot}\mathbf{e}_i^3 f_i^{eq}+\xi\Delta x^3\nabla \cdot \mathbf{e}_i \sum_{l=1}^q\Big(\mathbf{\Lambda}^{-1}-\frac{\mathbf{I}}{2}\Big)_{il}R_l
 	-\xi\Delta x^2\sum_{l=1}^q\Big(\mathbf{I}-\frac{\mathbf{\Lambda}}{2}\Big)_{il}R_l \Bigg)\notag\\
 	&
 	+\frac{\Delta x^2}{2}\mathbf{e}_k^2\overset{2}{\cdot}\nabla^2 \sum_{i=1}^q\Big( \mathbf{I}-\mathbf{\Lambda}^{-1}\Big)_{ki}\bigg(\xi\Delta x^2\partial_tf_i^{eq}+\Delta x\mathbf{e}_k\cdot\nabla f_i^{eq}-\Delta x^2\nabla^2\overset{2}{\cdot}\mathbf{e}_i\sum_{l=1}^q\mathbf{\Lambda}^{-1}_{il}\mathbf{e}_l f_l^{eq}\notag\\
 	&\qquad\qquad\qquad\qquad\qquad\qquad+\frac{\Delta x^2}{2}\mathbf{e}_i^2\overset{2}{\cdot}\nabla^2f_i^{eq}-\xi\Delta x^2\sum_{l=1}^q\Big(\mathbf{I}-\frac{\mathbf{\Lambda}}{2}\Big)_{il}R_l\bigg)\notag\\
 	&
 	-\frac{\Delta x^4}{6}\nabla^4\overset{4}{\cdot}\mathbf{e}_k^3 \sum_{i=1}^q\Big( \mathbf{I}-\mathbf{\Lambda}^{-1}\Big)_{ki}
 	\mathbf{e}_if_i^{eq} \notag\\
 	&
 	+\xi\Delta x^2\sum_{i=1}^q\Big(\mathbf{I}-\frac{\mathbf{\Lambda}}{2}\Big)_{ki}R_i
 	-\xi\Delta x^3\mathbf{e}_k\cdot\nabla\sum_{i=1}^q\Big(\mathbf{I}-\frac{\mathbf{\Lambda}}{2}\Big)_{ki}R_i
 	+\frac{\xi\Delta x^4}{2}\nabla^2\overset{2}{\cdot}\mathbf{e}_k^2\sum_{i=1}^q\Big(\mathbf{I}-\frac{\mathbf{\Lambda}}{2}\Big)_{ki}R_i .
 \end{align}  
 And  the unknown term $\sum_{k=1}^q\mathbf{e}_{ki}f_k^{ne}$ can be determined by taking the first-order moment  of the above equation [see (\ref{fneq-5rd})]. 
 \section{The Matlab symbolic computation code}\label{code}
  We  here provide a simple code to determine the  relaxation parameters and weight coefficients of the fourth order MRT-LB model (\ref{lb}) for the diagonal-anisotropic diffusion equation (\ref{DE}). To use the code for a high-dimensional case, one only needs to run the following code for each group consisting of two different directions, and the relaxation parameters and weight coefficients will be determined by the code output.
 
 \textbf{function} \quad [w$\_$xi,w$\_$xj,s$\_$xi,s$\_$xj,s$\_$xixj,s$\_$xi2xj,s$\_$xixj2] \quad ...

\vspace{0.2em} \qquad\qquad\qquad  $=$  Fun$\_$LB$\_$relax$\_$para(ww,d,s2i,eta,dt,dx,kappa$\_$xi,kappa$\_$xj) 

\vspace{0.2em} 
syms \quad sij\quad si \quad sj \quad wi\quad wj\quad real;

\vspace{0.2em} 
s$\_$2j  $=$  s$\_$2i; \quad si2j  $=$  sj;\quad sij2  $=$  si;

\vspace{0.2em} 
epsilon$\_$i  $=$  kappa$\_$xi*dt/dx$^2$ ;

\vspace{0.2em} 
epsilon$\_$j  $=$  kappa$\_$xj*dt/dx$^2$ ; 

\vspace{0.2em}     
eqi  $=$    epsilon$\_$i  $-$  (1/si $-$ 1/2)*(2*wi $+$ 4*(d $-$ 1)*ww); 

\vspace{0.2em} 
eqj  $=$   epsilon$\_$j  $-$  (1/sj $-$ 1/2)*(2*wj $+$ 4*(d $-$ 1)*ww);

\vspace{0.2em} 
eqii  $=$   ((1/2*epsilon$\_$i*(1/si $-$ 1/2)- \quad... 

\vspace{0.2em} 
\qquad \qquad   	( 
1/3/si $-$ 7/24 $+$ (1/si/s2i $-$ 1/2/si $-$ 1/2/s2i)*(1/si $-$ 1) $-$ 1/2*(1 $-$ 1/s2i)*(1/si $-$ 1/2) \quad... 

\vspace{0.2em} 
\qquad\qquad
$+$ ((1/si $-$ 1)*(1 $-$ 1/s2i $-$ 1/si) $+$ 1/2 $-$ 1/2/s2i)*epsilon$\_$i))/(1/si $-$ 1/2)); 

\vspace{0.2em} 
eqjj  $=$   ((1/2*epsilon$\_$j*(1/sj $-$ 1/2) $-$  \quad... 

\vspace{0.2em} 
\qquad\qquad  	(
1/3/sj $-$ 7/24 $+$ (1/sj/s2j $-$ 1/2/sj $-$ 1/2/s2j)*(1/sj $-$ 1) $-$ 1/2*(1 $-$ 1/s2j)*(1/sj $-$ 1/2) \quad... 

\vspace{0.2em} 
\qquad\qquad    	 $+$ ((1/sj $-$ 1)*(1 $-$ 1/s2j $-$ 1/sj) $+$ 1/2 $-$ 1/2/s2j)*epsilon$\_$j  	))/(1/sj $-$ 1/2)); 

\vspace{0.2em} 
eqij  $=$   (epsilon$\_$i*epsilon$\_$j $-$  \quad...

\vspace{0.2em} 
\qquad\qquad  	(
4*( $-$ 7/24 $+$ 1/6/si $+$ (1/si2j/s2i $-$ 1/2/si2j $-$ 1/2/s2i)*(1/si $-$ 1)\quad ... 

\vspace{0.2em} 
\qquad\qquad $-$ 1/2*(1 $-$ 1/s2i)*(1/si2j $-$ 1/2) $+$ 1/6/si2j)*ww \quad... 

\vspace{0.2em} 
\qquad\qquad
$+$ 4*( $-$ 7/24 $+$ 1/6/sj $+$ (1/sij2/s2j $-$ 1/2/sij2 $-$ 1/2/s2j)*(1/sj $-$ 1)\quad ... 

\vspace{0.2em} 
\qquad\qquad  $-$ 1/2*(1 $-$ 1/s2j)*(1/sij2 $-$ 1/2) $+$ 1/6/sij2)*ww\quad ... 

\vspace{0.2em} 
\qquad\qquad	 $+$ 4*( $-$ 7/24 $+$ 1/6/si $+$ (1/si2j/sij $-$ 1/2/si2j $-$ 1/2/sij)*(1/si $-$ 1)\quad ... 

\vspace{0.2em} 
\qquad\qquad $-$ 1/2*(1 $-$ 1/sij)*(1/si2j $-$ 1/2) $+$ 1/6/si2j)*ww \quad ... 

\vspace{0.2em} 
\qquad\qquad
$+$ 4*( $-$ 7/24 $+$ 1/6/sj $+$ (1/sij2/sij $-$ 1/2/sij2 $-$ 1/2/sij)*(1/sj $-$ 1)\quad ... 

\vspace{0.2em} 
\qquad\qquad $-$ 1/2*(1 $-$ 1/sij)*(1/sij2 $-$ 1/2) $+$ 1/6/sij2)*ww\quad ... 

\vspace{0.2em} 
\qquad\qquad
$+$ 4*( $-$ 7/24 $+$ 1/6/si $+$ (1/sij2/sij $-$ 1/2/sij2 $-$ 1/2/sij)*(1/si $-$ 1)\quad ... 

\vspace{0.2em} 
\qquad\qquad $-$ 1/2*(1 $-$ 1/sij)*(1/sij2 $-$ 1/2) $+$ 1/6/sij2)*ww \quad...

\vspace{0.2em} 
\qquad\qquad
$+$ 4*( $-$ 7/24 $+$ 1/6/sj $+$ (1/si2j/sij $-$ 1/2/si2j $-$ 1/2/sij)*(1/sj $-$ 1)\quad ... 

\vspace{0.2em} 
\qquad\qquad $-$ 1/2*(1 $-$ 1/sij)*(1/si2j $-$ 1/2) $+$ 1/6/si2j)*ww \quad... 

\vspace{0.2em} 
\qquad\qquad
$+$ ((1/si $-$ 1)*(1 $-$ 1/s2i $-$ 1/si) $+$ 1/2 $-$ 1/2/s2i)*epsilon$\_$j*(2*wi $+$ 4*(d $-$ 1)*ww) \quad... 

\vspace{0.2em} 
\qquad\qquad
$+$ ((1/sj $-$ 1)*(1 $-$ 1/s2j $-$ 1/sj) $+$ 1/2 $-$ 1/2/s2j)*epsilon$\_$i*(2*wj $+$ 4*(d $-$ 1)*ww)  	));  

\vspace{0.2em} 
h $=$  solve(eqi,eqj,eqii,eqjj,eqij,wi,wj,si,sj,sij);

\vspace{0.2em} 
Wi $=$ double(h.wi);\quad Wj $=$ double(h.wj);

\vspace{0.2em} 
Si $=$ double(h.si);\quad Sj $=$ double(h.sj);\quad Sij $=$ double(h.sij); 

\vspace{0.2em} 
w$\_$xi=[ ];\quad w$\_$xj=[ ];\quad s$\_$xi$\_$=[ ];\quad s$\_$xj$\_$=[ ]; \quad s$\_$xixj=[ ];

\vspace{0.2em} 
\textbf{for} i=1:length(Wi)

\vspace{0.2em} 
\qquad \textbf{if} Wi(i)$>0$ \&\& Wi(i)$<1$ \&\& Wj(i)$>0$ \&\& Wj(i)$<1$  \&\& \quad ...

\vspace{0.2em} 
\qquad \qquad 
Si(i)$>0$ \&\& Si(i)$<2$ \&\& Sj(i)$>0$ \&\& Sj(i)$<2$  \&\& Sij(i)$>0$ \&\& Sij(i)$<2$  

\vspace{0.2em} 
\qquad \qquad w$\_$xi $=$ [w$\_$xi;  Wi(i)];\quad  w$\_$xj $=$ [w$\_$xj; Wj(i)];

\vspace{0.2em} 
\qquad \qquad s$\_$xi$\_$ $=$ [s$\_$xi$\_$; Si(i)];\quad  s$\_$xj$\_$  $=$ [s$\_$xj$\_$; Sj(i)];\quad s$\_$xixj $=$ [s$\_$xixj; Sij(i)];

\vspace{0.2em} 
\qquad \textbf{break;}

\vspace{0.2em} 
\textbf{end}

\vspace{0.2em} 
s$\_$xi  $=$  2*eta*dt./(dt*eta $-$ sqrt((s$\_$xi$\_$ $-$ 4*dt*eta $+$ dt$^2$*s$\_$xi$\_$*eta$^2$ $+$ 2*dt*s$\_$xi$\_$*eta)./s$\_$xi$\_$) $+$ 1);

\vspace{0.2em} 
s$\_$xj  $=$  2*eta*dt./(dt*eta $-$ sqrt((s$\_$xj$\_$ $-$ 4*dt*eta $+$ dt$^2$*s$\_$xj$\_$*eta$^2$ $+$ 2*dt*s$\_$xj$\_$*eta)./s$\_$xj$\_$) $+$ 1);

\vspace{0.2em} 
s$\_$xi2xj  $=$  s$\_$xj; \quad s$\_$xixj2 $=$ s$\_$xi;
 
	\section*{Acknowledgments}
	The computation is completed in the HPC Platform of Huazhong University of Science and Technology. This work was financially supported by the National Natural Science Foundation of China (Grants No. 12072127 and No. 123B2018), the Interdisciplinary Research Program of HUST (Grants No. 2023JCJY002 and No. 2024JCYJ001), and the Fundamental Research Funds for the Central Universities of HUST  (Grants No. 2024JYCXJJ016 and No. YCJJ20241101).
	
 \bibliographystyle{elsarticle-harv} 
 \bibliography{reference}
 
\end{document}